\documentclass[11pt]{article}
\usepackage{amssymb,amsmath,latexsym}
\usepackage{epsfig}

\newtheorem{thm}{Theorem}[section]

\newtheorem{lemma}[thm]{Lemma}
\newtheorem{prop}[thm]{Proposition}
\newtheorem{defn}[thm]{Definition}

\newtheorem{remark}[thm]{Remark}
\numberwithin{equation}{section}

\def\pf{{\medskip\noindent {\bf Proof. }}}
\def\qed{{\hfill $\Box$ \bigskip}}

  \def\sC {{\cal C}}
 \def\sE {{\cal E}} \def\sF {{\cal F}}
  
  \def\sL {{\cal L}}

 \def\bK {{\mathbb K}} 
  
  \def\bR {{\mathbb R}}

\def\wt{\widetilde}

\def\E{{\mathbb E}}
\def\P{{\mathbb P}}

\def\Exp{{\rm Exp}}

\def\nn{\nonumber}

\def\bea{\begin{align*}}
\def\eea{\end{align*}}
\def\bee{\begin{equation}}
\def\eee{\end{equation}}

\def\eps{\varepsilon}

\makeatletter
\@addtoreset{equation}{section}

\makeatother

\begin{document}
\allowdisplaybreaks
\bibliographystyle{plain}

\title{\Large \bf Stability of Dirichlet heat kernel estimates for non-local operators under Feynman-Kac perturbation}

\author{{\bf Zhen-Qing Chen}\thanks{Research partially supported
by NSF Grants DMS-0906743 and DMR-1035196.}, \quad {\bf Panki Kim}\thanks{This research was supported by Basic Science Research Program through the National Research Foundation of Korea(NRF) funded by the Ministry of Education, Science and Technology(0409-20110087).} \quad and \quad {\bf Renming
Song}\thanks{Research supported in part by a grant from the Simons
Foundation (208236).} }
\date{(December 11, 2011)}
\maketitle

\begin{abstract}
In this paper we show that
Dirichlet heat kernel estimates for a class
of (not necessarily symmetric) Markov processes are stable under non-local
Feynman-Kac perturbations. This class of processes includes, among others,
(reflected) symmetric stable-like processes on closed $d$-sets in $\bR^d$,
killed symmetric stable processes, censored stable processes in $C^{1, 1}$ open sets
as well as stable processes with drifts in bounded
$C^{1, 1}$ open sets.
\end{abstract}

	

\bigskip
\section{Introduction}\label{s:1}

Recently, sharp two-sided Dirichlet heat kernel estimates have been
obtained for several classes of discontinuous processes
(or non-local operators), including
symmetric stable processes \cite{CKS}, censored stable processes
\cite{CKS1}, relativistic stable processes \cite{CKS2}, and stable processes with drifts \cite{CKS4}.
Although the proofs in these papers share a general
road map, there are many distinct difficulties
and the actual arguments are specific to the underlying processes.
The main purpose of this paper is to establish a stability result
for the sharp Dirichlet heat kernel estimates of a family of discontinuous
processes under non-local Feynman-Kac perturbations.
Here for a discontinuous Hunt process $X$, a non-local Feynman-Kac
transform is given by
$$ T_t f (x) =\E_x \Big[
  \exp \Big(A_t +\sum_{s\leq t} F(X_{s-}, X_s) \Big)f(X_t)\Big],
$$
where $A$ is a continuous additive functional of $X$ having finite
variations each compact time interval
and $F(x, y)$ is a measurable function that
vanishes along the diagonal.
The approach of this paper is quite robust that it applies to a
class of not necessarily symmetric Markov processes which includes
all the four families of processes mentioned above
in bounded $C^{1,1}$ open sets.

Transformation by multiplicative functionals is
one of the most important transforms for Markov processes
(see, for example, \cite{Chu, sharpe}).
Non-local Feynman-Kac transforms are particular cases.
They  play an important role in the probabilistic
as well as analytic aspect of potential theory, and also in mathematical physics.
For example, it is shown in \cite{CS03b} that relativistic
stable processes can be obtained from the symmetric
$\alpha$-stable processes through Feynman-Kac transformations.
We refer the reader to \cite{ChZ, Sim} for nice accounts
on Feynman-Kac semigroups of Brownian motion.
In particular, it is shown in
\cite{BM, Sim}
 that under a certain
Kato class condition, the integral kernel (called the heat kernel)
of the Feynman-Kac semigroup
of Brownian motion admits two-sided Gaussian bound estimates.
In \cite{R}, sharp two-sided estimates on the  densities of (local)
Feynman-Kac semigroups of killed Brownain motions in $C^{1, 1}$
domains were established.
Non-local Feynman-Kac semigroups for
symmetric stable processes and their associated quadratic forms
were studied in \cite{S1, S2}.
By combining some ideas from \cite{Zq2} with
results from \cite{CK}, it was proved in \cite{S} that,
under a certain Kato class condition, the heat kernel of
the non-local Feynman-Kac semigroup of a symmetric
stable-like process $X$ on $\bR^d$ is comparable to that of $X$.
The symmetry condition on $F(x, y)$ plays an essential role
in the argument of \cite{S}.
The nonsymmetric pure jump case for stable-like processes
is dealt with in \cite{W}.
For recent development in the study of non-local Feynman-Kac transforms
for general symmetric Markov processes,
we refer the reader to \cite{CFKZ1, CFKZ2} and the references therein.
We also mention that the stability of Martin boundary under non-local Feynman-Kac perturbation is addressed in \cite{CKi}.
To the best of the authors knowledge, Dirichlet heat kernel estimates for (either local or non-local) Feynman-Kac semigroups
of discontinuous processes is studied here for the first time.
The main challenge in studying Dirichlet heat kernel estimates
of Feynman-Kac semigroups is to get exact boundary decay
behavior of the heat kernels.
While our main interest is in the Dirichlet heat kernel estimates
for Feynman-Kac semigroups, our theorem also covers the whole space
 case as well as ``reflected" stable-like processes on subsets
of $\bR^d$.
In particular, our result recovers and extends
 the main results of \cite{S, W}
where $D=\bR^d$.
Even in the whole space case, our approach is different from
those in \cite{S, W}.

\subsection{Setup and main result}

In this paper we
always assume that $\alpha\in (0, 2)$, $d \ge 1$,
$D$ is
a Borel
set in $\bR^d$.
For any $x\in D$, $\delta_D(x)$ denotes the
Euclidean distance between $x$ and $D^c$.
We use ``$:=$" to denote a definition, which is
read as ``is defined to be".
For $a, b\in \bR$,
$a\wedge b:=\min \{a, b\}$ and $a\vee b:=\max \{a, b\}$.
 The Euclidean distance between $x$ and $y$ is
denoted as $|x-y|$.

For $\gamma \geq 0$, let
$$
\psi_\gamma(t, x, y):= \left(1\wedge
\frac{\delta_D(x)}{t^{1/\alpha}}\right)^\gamma \left(1\wedge
\frac{\delta_D(y)}{t^{1/\alpha}}\right)^\gamma,
\qquad t>0, \, x, y \in D.
$$
Throughout this paper, $X$ is a Hunt process
on $D$
with transition semigroup
$\{ P_t: t\ge 0\}$ that  admits a jointly continuous
transition density $p_D(t, x, y)$ with respect to
the Lebesgue measure and that there exist
$C_0>1$ and $\gamma\in [0, \alpha\wedge d)$ such that
\begin{equation}\label{e:hke}
C_0^{-1} \psi_\gamma(t, x, y)q(t,x,y) \le p_D(t, x, y) \le C_0
\psi_\gamma(t, x, y)q(t,x,y)
\end{equation}
for all $(t, x, y)\in (0, 1]\times D\times D$,
where
\begin{equation}\label{e:defq}
q(t,x, y)
:=t^{-d/\alpha}\wedge\frac{t}{|x-y|^{d+\alpha}}.
\end{equation}
It is easy to see that under this assumption, $X$ is a Feller process
satisfying the strong Feller property.
It is easy to see that, by increasing the value of $C_0$ if necessary,
\begin{equation}\label{e:hkein1}
C_0^{-1} \le \int_{\bR^d} q(t,x,y)dy \le C_0
\quad \text{ for all }(t, x, y)\in (0, \infty)\times \bR^d.
\end{equation}
Thus
\begin{equation}\label{e:hkein2}
 \int_{D} p_D(t, x, y) dy \le C_0^2 \left(1\wedge
\frac{\delta_D
(x)}{t^{1/\alpha}}\right)^\gamma
 \quad \text{ for all }(t, x)\in (0, 1]\times D.
\end{equation}

Note that $X$ is not necessarily symmetric.
We further assume that $X$ has a L\'evy system $(N, t)$ where
$N=N(x, dy)$ is a kernel given by
$$
N(x, dy)=\frac{c(x, y)}{|x-y|^{d+\alpha}}dy,
$$
with $c(x, y)$ a measurable function that is bounded between
two positive constants on $D\times D$.
That is, for
any  $x\in D$, any stopping time $T$ (with respect to the filtration of $X$)
and any non-negative measurable function $f$ on $D\times D$
with $f(y, y)=0$ for all $y\in D$ that is extended to be zero off $D\times D$,
\begin{equation}\label{e:levy}
\E_x
\left[\sum_{s\le T} f(X_{s-}, X_s) \right]= \E_x \left[
\int_0^T \left( \int_{D} f(X_s, y) \frac{c(X_s, y)}{|X_s-y|^{d+\alpha}} dy \right) ds \right].
\end{equation}
By increasing the value of $C_0$ if necessary, we may and do assume
that
\begin{equation}\label{e:1.5a}
 1/C_0 \leq c(x, y) \leq C_0 \qquad \hbox{for } x, y\in D.
\end{equation}

Recall that an open set $D$ in $\bR^d$ (when $d\ge 2$) is said to be
a $C^{1,1}$ open set if there exist a localization radius $ r_0>0 $
and a constant $\Lambda_0>0$ such that for every $z\in\partial D$,
there exist a $C^{1,1}$-function $\phi=\phi_z: \bR^{d-1}\to \bR$
satisfying $\phi (0)=0$, $\nabla\phi (0)=(0, \dots, 0)$, $\| \nabla
\phi  \|_\infty \leq \Lambda_0$,
$| \nabla \phi (x)-\nabla \phi (w)|
\leq \Lambda_0 |x- w|$, and an orthonormal coordinate system $y=(y_1,
\cdots, y_{d-1}, y_d):=(\wt y, \, y_d)$ such  that $ B(z, r_0 )\cap
D=B(z, r_0 )\cap \{ y: y_d > \phi (\wt y) \}$.
We call the pair $(r_0, \Lambda_0)$ the
characteristics of the $C^{1,1}$ open set $D$.
By a $C^{1,1}$ open set
in $\bR$ we mean an open set which can be
 expressed
 as the union of disjoint intervals so that the minimum of
the lengths of all these intervals is positive and the minimum of
the distances between these intervals is positive.

\medskip

It follows from
\cite{CKS, CKS1, CKS4, CK} that the following are true:
\begin{description}
\item{(i)} the (reflected) symmetric stable-like process  on
any closed $d$-subset $D$ in $\bR^d$ (see Subsection \ref{S:4.1} for the definition of   $d$-set)
satisfies the conditions \eqref{e:hke} and \eqref{e:levy} with $\gamma=0$ and  $c(x, y)$
a symmetric measurable function that is bounded between two positive constants;

\item{(ii)}
 the killed symmetric $\alpha$-stable process
on a $C^{1, 1}$ open set $D$ satisfies the conditions \eqref{e:hke} and \eqref{e:levy} with $\gamma=\alpha/2$ and $c(x, y)=c$;

\item{
(iii)} when $d \ge 2$ and $\alpha \in (1,2)$, the killed symmetric $\alpha$-stable process
with drift in a bounded $C^{1, 1}$ open set $D$
satisfies the conditions \eqref{e:hke} and \eqref{e:levy} with
$\gamma=\alpha/2$ and $c(x, y)=c$; and

\item{(iv)}
when $\alpha \in (1,2)$,
the censored $\alpha$-stable process in a $C^{1, 1}$
open set $D$ satisfies the conditions \eqref{e:hke} and \eqref{e:levy}
with $\gamma=\alpha-1$ and $c(x, y)=c$.
\end{description}

By a signed measure $\mu$ we mean in this paper the difference of two
nonnegative $\sigma$-finite measures $\mu_1$ and $\mu_2$ in $D$.
We point out that
 $\mu=\mu_1-\mu_2$ may not be a signed measure in $D$
in the usual sense as both $\mu_1 (D)$ and $\mu_2(D)$ may be infinite.
However, there is an increasing sequence of subsets $\{F_k, k\geq 1\}$
whose union is $D$ so that $\mu_1(F_k)+\mu_2(F_k)<\infty$ for every
$k\geq 1$. So when restricted to each $F_k$, $\mu$ is a finite signed
measure. Consequently, the positive and negative parts of $\mu$
are well defined on each $F_k$ and hence on $D$, which will be
denoted as    $\mu^+$ and $\mu^-$, respectively.
We use $|\mu|=\mu^++\mu^-$ to denote the total variation measure of $\mu$.
Taking such an extended view of signed measures is desirable
when one studies the correspondence between signed measures and
continuous functions of finite variations or the correspondence
between signed smooth measures and continuous additive functionals
of finite variations for a Hunt process.
For a signed measure $\mu$ on $D$ and $t>0$,
we define
$$
N^{\alpha, \gamma}_{\mu}(t)=\sup_{x\in
D}\int^t_0\int_{D}\left(1\wedge
\frac{\delta_D(y)}{s^{1/\alpha}}\right)^\gamma q(s, x,
y)|\mu|(dy)ds.
$$

\begin{defn}\label{d:kc}
A signed measure $\mu$ on $D$ is said to be in the Kato
class ${\bf K}_{\alpha, \gamma}$ if \, $\lim_{t\downarrow 0}N^{\alpha,
\gamma}_{\mu}(t)=0$.
\end{defn}

Note that if $N^{\alpha, \gamma}_\mu(t)<\infty$ for some $t>0$,
then $|\mu|$ is a Radon measure on $D$.
We say that a measurable function
$
g$ belongs to the Kato class ${\bf K}_{\alpha, \gamma}$ if $
g(x)dx
\in {\bf K}_{\alpha, \gamma}$
and we denote $N^{\alpha, \gamma}_{
g (x) dx}$ by $N^{\alpha, \gamma}_
g$.
It is well known that any $\mu\in {\bf K}_{\alpha, \gamma}$ is a
smooth measure in the sense of \cite{FG}.
Moreover,
using the fact that $X$ has a transition density function under
each $\P_x$, one can show that the continuous
additive functional $A^{\mu}_t$ of $X$ with Revuz measure $\mu \in {\bf K}_{\alpha, \gamma}$ can
be defined without exceptional set, see \cite[pp. 236--237]{FOT} for details.
Concrete conditions for $\mu\in {\bf K}_{\alpha, \gamma}$ are given
in Proposition \ref{P:4.1}.

For any
measurable
 function $F$ on $D\times D$ vanishing on the
diagonal, we define
$$
N^{\alpha, \gamma}_F(t) := \sup_{y\in D}
\int_0^t \int_{
D\times D}
\Big(1\wedge \frac{\delta_D(z)}{s^{1/\alpha}}\Big)^\gamma \,
q(s, y, z ) \,
\Big(1 + \frac{|z-w| \wedge t^{1/\alpha}}{|y-z|} \Big)^\gamma
\, \frac{|F|(z,w) + |F|(w,z)}{|z-w|^{d+\alpha}} dw dz ds.
$$

\begin{defn}\label{d:nlkc}
Suppose that $F$ is a
 measurable
 function  on $D\times D$ vanishing
on the diagonal. We say that $F$ belongs to the Kato class ${\bf
J}_{\alpha, \gamma}$
 if $F$ is bounded and
$\lim_{t\downarrow 0}N^{\alpha, \gamma}_{F}(t)=0$.
\end{defn}

It follows immediately from the two definitions above that if $F\in
{\bf J}_{\alpha, \gamma}$,
then the function
$$
z\mapsto \int_D\frac{|F|(z,w) + |F|(w,z)}{|z-w|^{d+\alpha}} dw
$$
belongs to ${\bf K}_{\alpha, \gamma}$.
See Proposition \ref{p:newqw1} for a sufficient condition for $F\in
{\bf J}_{\alpha, \gamma}$.

It is easy to check that if $F$ and $G$ belong to ${\bf J}_{\alpha,
\gamma}$ and $c$ is a constant, then the functions $cF, e^{F}-1,
F+G$ and $FG$ all belong to ${\bf J}_{\alpha, \gamma}$.
Throughout this paper, we will use the following notation: For any
given measurable function $F$
on $D\times D$,  $F_1(x, y)$
denotes the function $e^{F(x, y)}-1$.

For any signed measure $\mu$ on $D$ and any
measurable function $F$ on
$D\times D$ vanishing on the diagonal, we define
$$
N^{\alpha, \gamma}_{\mu, F}(t) :=N^{\alpha, \gamma}_{\mu}(t)+
N^{\alpha, \gamma}_F(t).
$$
When $\mu\in{\bf K}_{\alpha, \gamma}$ and
$F$ is
a measurable  function with $F_1\in {\bf J}_{\alpha, \gamma}$,
we put
$$
A^{\mu, F}_t
=A^\mu_t+\sum_{0<s\le t}F(X_{s-}, X_s).
$$
For any nonnegative Borel function $f$ on $D$, we define
$$
T^{\mu, F}_t f(x)
=\E_x\left[\exp(A^{\mu, F}_t) f(X_t)\right], \quad t\ge 0,
x\in D.
$$
Then $(T^{\mu, F}_t: t\ge 0)$ is called the Feynman-Kac semigroup of $X$ corresponding to $\mu$ and $F$.
The main purpose of this paper is to establish the following
result. Recall that $\gamma\ge 0$ and $C_0\geq 1$ are the constants
in \eqref{e:hke} and \eqref{e:1.5a}.
For any bounded function $F$ on $D\times D$, we use $\|F\|_\infty$
to denote $\|F\|_{L^\infty (D \times D)}$.

\begin{thm}\label{t:main}
 Let $d\geq 1$, $\alpha \in (0, 2)$ and
 $\gamma \in [0, \alpha \wedge d)$.
 Suppose $X$ is a Hunt process in
a Borel set $D\subset \bR^d$
with a jointly continuous
transition density $p_D(t, x, y)$ satisfying \eqref{e:hke}, \eqref{e:levy} and \eqref{e:1.5a}.
If $\mu$ is a signed measure in ${\bf K}_{\alpha, \gamma}$ and $F$ is
a measurable  function so that $F_1:=e^F-1\in {\bf J}_{\alpha, \gamma}$, then the
non-local Feynman-Kac semigroup $(T^{\mu, F}_t: t\ge 0)$
has a continuous
density $q_D(t, x, y)$, and for any $T>0$,
there exists a constant
$C=C(d, \alpha, \gamma, C_0,
N^{\alpha, \gamma}_{\mu, F_1}, \|F_1\|_\infty,
T) >0$ such that for all
$(t, x, y)\in (0, T]\times D\times D$,
$$
 q_D(t, x, y) \le C \psi_\gamma(t, x, y)q(t,x,y).
$$
If $\mu \in {\bf K}_{\alpha, \gamma}$
and $F\in {\bf J}_{\alpha, \gamma}$, then
there exists a constant $\wt C= \wt C(d, \alpha, \gamma, C_0,
N^{\alpha, \gamma}_{\mu, F},
\|F\|_\infty, T) >1$ such that for all
$(t, x, y)\in (0, T]\times D\times D$,
$$
 \wt C^{-1} \psi_\gamma(t, x, y)q(t,x,y) \leq
 q_D(t, x, y) \le \wt C \psi_\gamma(t, x, y)q(t,x,y).
$$
\end{thm}

Here and in the sequel, the dependence of the constant $C$ on $N^{\alpha, \gamma}_{\mu, F_1}$ and $\|F_1\|_\infty$ means that the value of the constant
$C$ depends only on a specific upper bound for the rate
of the function $N^{\alpha, \gamma}_{\mu, F_1}(t)$
going
to zero as $t\to 0$ and on a specific upper bound for $\|F_1\|_\infty$.
When $D=\bR^d$ and $\gamma=0$,
Theorem \ref{t:main} in particular recovers and extends
the main results of \cite{S, W}.

\subsection{Approach}

To explain our approach, we first recall the definition of the Stieltjes exponential.
If $K_t$ is a right continuous function with left limits on $\bR_+$ with
$K_0=1$ and $\Delta K_t:=K_t-K_{t-}>-1$ for every $t>0$, and if
$K_t$ is of finite variation on each compact time interval, then the
Stieltjes exponential $\Exp (K)_t$ of $K_t$ is the unique
solution $Z_t$ of
$$
 Z_t=1+ \int_{(0, t]} Z_{s-} dK_s, \quad t> 0.
$$
By \cite[(A4.17)]{sharpe},
\begin{equation}\label{e:stie0}
 \Exp (K)_t = e^{K^c_t} \prod_{0<s\leq t} (1+\Delta K_s),
\end{equation}
where $K_t^c$ denotes
the continuous part of $K_t$.
Clearly $\exp (K_t)\geq \Exp (K)_t$ with the equality
holds if and only if $K_t$ is continuous.
The reason of $\Exp (K)_t$ being called the {\it Stieltjes}
exponential of $K_t$ is that  by  \cite{DD}
we have
\begin{equation}\label{e:stie1}
 \Exp (K)_t =1+ \sum_{n=1}^\infty
\int_{(0, t]} dK_{t_n} \int_{(0, t_n]} dK_{t_{n-1}}
\cdots \int_{(0, t_2]} dK_{t_1}.
\end{equation}
The advantage of using the Stieltjes exponential $\Exp (K)_t$
over the usual exponential $\exp (K_t)$  is the identity
\eqref{e:stie1}, which allows one to apply
the Markov property of $X$.

Recall that $F_1(x, y)=e^{F(x, y)}-1$.
In view of \eqref{e:stie0}, we can express $\exp (A^{\mu, F}_t)$
in terms of the Stieltjes exponential:
$$
\exp (A^{\mu, F}_t )  = \Exp \Big( A^\mu+ \sum_{s\leq \cdot }
F_1 (X_{s-}, X_s)  \Big)_t \qquad \hbox{for } t\geq 0.
$$
Applying \eqref{e:stie1} with $K_t:=A^\mu_t+\sum_{s\leq t}
F_1 (X_{s-}, X_s)$ and using the Markov property of $X$,
 we have for any  bounded $f\geq 0$ on $D$,
\begin{eqnarray}
T^{\mu, F}_t f (x)
&=& \E_x\left[\exp (A^{\mu, F}_t )f(X_t)\right]
=\E_x \left[f(X_t) \,\Exp \Big( A^\mu + \sum_{s\leq \cdot }
F_1 (X_{s-}, X_s) \Big)_t  \right] \nonumber   \\
&=& P_t f(x) + \E_x \left[ f(X_t)\sum_{n=1}^\infty
\int_{(0, t]} dK_{t_n} \int_{(0, t_n]} dK_{t_{n-1}}
\cdots \int_{(0, t_2]} dK_{t_1}\right] . \label{e:newq3}
\end{eqnarray}

It can be shown that, for $\mu\in{\bf
K}_{\alpha, \gamma}$ and
measurable function $F$ with $F_1
\in {\bf J}_{\alpha, \gamma}$,
there is some constant $T_0>0$ so that we can
change the order of the expectation and the infinite sum when $t\leq T_0$.
Hence we have for every $t\leq T_0$,
\begin{eqnarray}\label{e:1.11}
T^{\mu, F}_t f(x) &=&  P_t f(x) + \sum_{n=1}^\infty \E_x \left[ f(X_t)
\int_{(0, t]} dK_{t_n} \int_{(0, t_n]} dK_{t_{n-1}}
\cdots \int_{(0, t_2]} dK_{t_1}\right] \nonumber \\
&=&P_t f(x) + \sum_{n=1}^\infty \E_x \left[
\int_{(0, t]} P_{t-t_n} f(X_{t_n})dK_{t_n} \int_{(0, t_n]} dK_{t_{n-1}}
\cdots \int_{(0, t_2]} dK_{t_1} \right].
\end{eqnarray}
Note that by \eqref{e:levy},  for any bounded function $g$,
\begin{align}\label{e:1.12}
& \E_x \left[ \int_{(0, s]} g(X_{r}) dK_r \right]
=\E_x \left[ \int_{(0, s]} g(X_{r}) dA^\mu_r
+ \sum_{r\leq s}
g(X_{r}) F_1(X_{r-}, X_r)\right] \nonumber \\
=& \int_0^s \int_D p_D(r, x, y) g(y) \mu  (dy) dr
+ \E_x \left[
\int_0^s \left( \int_{D}  F_1(X_r, y)g(y)  \frac{c(X_r, y)}{|X_r-y|^{d+\alpha}} dy \right) dr \right] \nonumber \\
=&\int_0^s \int_D p_D(r, x, y) g(y) \mu  (dy) dr
+
\int_0^s \int_D p_D(r, x, z) \left( \int_{D} F_1(z, y)  g(y)
 \frac{c(z, y)}{|y-z|^{d+\alpha}} dy \right) dz dr .\end{align}
This together with \eqref{e:newq3} motives us to define
 $p^0 (t, x, y):=p_D(t, x, y)$
and, for $k \ge 1$
\begin{eqnarray}
p^k (t, x, y)&=&\int^t_0 \left(\int_{D}p_D(s, x,
z)p^{k-1}(t-s, z, y)\mu(dz) \right)ds\nonumber \\
&&+ \int^t_0 \left(\int_{D\times D} p_D( s, x,
z)\frac{c(z, w) F_1(z, w)}{|z-w|^{d+\alpha}}p^{k-1}(t-s, w, y) dzdw
\right)ds.\label{e:i0}
\end{eqnarray}
One then concludes from \eqref{e:newq3} that
$$ T^{\mu, F}_t f(x) = \int_D q_D(t, x, y) f(y) dy,
$$
where
\begin{equation}\label{e:1.14}
q_D(t, x, y) := \sum_{k=0}^\infty p^k(t, x, y).
\end{equation}
We then proceed to establish the following key estimates:
there exist constants $T_1\in (0, T_0]$, $c>0$ and $0<\lambda <1$
such that
\begin{equation}\label{e:1.15}
 |p^k(t, x, y)| \leq (\lambda^k+ck \lambda^{k-1} )p_D(t, x, y)
\qquad \hbox{on } (0, T_1] \times D \times D
\hbox{ for every } k\geq 1.
\end{equation}
From this we can deduce  that  for every $t\in (0, T_1]$,
\begin{equation}\label{e:1.15a}
 q_D(t, x, y)=
\sum_{k=0}^\infty p^k(t, x, y)
\leq \left(\frac1{1-\lambda}+ \frac{c}{(1-\lambda)^2}\right)\, p_D(t, x, y),
\end{equation}
and, under the assumption $F\in {\bf J}_{\alpha, \gamma}$, that
$$
q_D(t, x, y)   \ge 2^{-2(\lambda+c)}p_D(t, x, y), $$
which establish Theorem \ref{t:main} for $t\leq T_1$.
The general case of $t\leq T$ follows from an application
 of the Chapman-Kolmogorov equation.

The key to establish the estimate \eqref{e:1.15}
are two integral forms of the 3P inequality given in Lemma
\ref{l:3p} and Theorem \ref{t:G3p} below.
For a killed Brownian motion in a smooth domain, the following form of 3P inequality is known (see \cite{KS1, R1}):
for any $0< c < a \wedge (b-a)$, there exists $M=M(a,b,c)>0$
such that for every $0<s<t$ and $x,y,z \in D$
\begin{equation}\label{e:ndf1}
\frac{p^W_a(t-s, x,
z) p^W_b(s, z,y)}{p^W_a(t, x, y)}  \le M \frac{\delta_D(z)}{\delta_D(x)}
p^W_c(t-s, x, z)+ M\frac{\delta_D(z)}{\delta_D
(y)} p^W_c(s, y, z)
\end{equation}
where $ p^W_c(t,x,y) := \psi_1(t, x, y) t^{-d/2} e^{-c|x-y|^2/t}. $
For symmetric $\alpha$-stable processes in $\bR^d$,
one has  the following form of 3P inequality (see \cite{BJ} and \eqref{e:ppp} below):
\begin{equation}\label{e:ndf2}
\frac{ q(s,x,z) q(t-s, z,y)}{ q(t,x,y)}\le
c\left(q(s,x,z)+ q(t-s, z,y) \right) \quad \text{for every }0<s<t \text{ and } x,y,z \in \bR^d .
\end{equation}
The above 3P type inequalities \eqref{e:ndf1} and \eqref{e:ndf2}
played essential roles in establishing the  heat kernel estimates in \cite{BJ, KS1, R1}.
It seems that, for the processes we are dealing with in this paper,
the above two types of 3P inequalities are not true in general.
Moreover, we need a 3P type estimate on  $p_D(t-s, x,
z) p_D(s, w,y)/p_D(t, x, y)$, where $z \not= w$.

\bigskip

The rest of the paper is organized as follows.
In Section
\ref{s:2}, we prove some key inequalities, including two forms of the 3P inequality.
The main estimates \eqref{e:1.15} and Theorem \ref{t:main} will be established in Section \ref{s:3}.
In the last section, we give some applications of our
main results.

In this paper, we will use capital letters $
\wt C, C, C_0,
C_1, C_2, \dots$ to denote constants in the statements of results,
and their values will be fixed. The lower case letters $c_1, c_2,
\dots$ will denote generic constants used in proofs, whose exact
values are not important and can change from one appearance to
another. The labeling of the lower case constants starts anew in
each proof.
For two positive functions $f$ and $g$,
we use the notation $f\asymp g$, which means that there are two positive constants
$c_1$ and $c_2$ whose
values depend only on $d, \alpha$ and $\gamma$ so that
$c_1 g\leq f \leq c_2 g$.

\section{3P inequalities}\label{s:2}

In this section we will establish some key inequalities which will
be essential in proving Theorem \ref{t:main}.
The main results of this section are Lemma
\ref{l:4-1}, Theorem \ref{t:key}, Lemma \ref{l:4-2}
 and Theorem \ref{t:G3p}.
Throughout this section, $D$ is a Borel set in $\bR^d$.

The following elementary facts will be used several times in this section.

\begin{lemma}\label{l:ineq2-1}
For any $s, t>0$ and $(y,z) \in D \times D$, we have
\begin{equation}\label{e:elineq1}
1\wedge \frac{\delta_D(z)}{t^{1/\alpha}} =
  \frac{\delta_D(y)}{t^{1/\alpha}}\,
\left(\frac{\delta_D(z)\wedge
t^{1/\alpha}}{\delta_D(y)}\right)
\end{equation}
and
\begin{equation}\label{e:elineq2}
\left(1\wedge \frac{\delta_D(y)}{s^{1/\alpha}}\right) \left(1\wedge
\frac{\delta_D(z)}{t^{1/\alpha}}\right) \le
2 \left(1+\frac{|y-z|}{s^{1/\alpha} +\delta_D (y)}\right)\left(1 \wedge
\frac{\delta_D(y)}{t^{1/\alpha}}\right) .
\end{equation}
\end{lemma}

\pf
The identity \eqref{e:elineq1} is clear, so we only need to prove \eqref{e:elineq2}.
Since $\delta_D(z) \le |y-z|+\delta_D(y)$, we see that
\begin{eqnarray*}
1\wedge \frac{\delta_D(z)}{t^{1/\alpha}}
\le  1\wedge
\left(\left(\frac{|y-z|+\delta_D(y)}{\delta_D(y)}\right)\frac{\delta_D(y)}{t^{1/\alpha}}\right)\le
 \left(1+\frac{|y-z|}{\delta_D(y)}\right)
\left(1\wedge \frac{\delta_D(y)}{t^{1/\alpha}}\right).
\end{eqnarray*}
Thus, applying the elementary inequality
\begin{equation}\label{e:ineq*}
\frac{a}{a +b} \le  1\wedge\frac{a}{b}
\le  \frac{2a }{a+ b}, \quad a, b>0
\end{equation}
we get
\begin{eqnarray*}
\left(1\wedge
\frac{\delta_D(y)}{s^{1/\alpha}}\right) \left(1\wedge
\frac{\delta_D(z)}{t^{1/\alpha}}\right) &\le&   \left(1\wedge
\frac{\delta_D(y)}{s^{1/\alpha}}\right)\left(1+
\frac{|y-z|}{\delta_D(y)}\right)
\left(1\wedge \frac{\delta_D(y)}{t^{1/\alpha}}\right)\\
&\le&2\left(1+\frac{|y-z|}{s^{1/\alpha}+\delta_D(y)}\right)
\left(1 \wedge \frac{\delta_D(y)}{t^{1/\alpha}}\right).
\end{eqnarray*}
 \qed

Using \eqref{e:defq} and
\eqref{e:ineq*}, we get that
\begin{equation}\label{e:ineq1}
\frac{t}{(t^{1/\alpha}+|x-y|)^{d+\alpha}}
\le q(t,x,y) \le 2^{d+\alpha}  \frac{t}{(t^{1/\alpha}+|x-y|)^{d+\alpha}}.
\end{equation}

\begin{lemma}\label{l:4-1}
 For any $\gamma\in [0, 2\alpha)$,
 there exists a constant
$C_1:=C_1(d, \alpha, \gamma)>1$
 such that for all $(t, y,z )\in (0,
\infty)\times D \times D$,
\begin{equation} \label{e:ineqmain1}
\left(1\wedge \frac{\delta_D(z)}{t^{1/\alpha}}\right)^\gamma
\int^{t/2}_0\psi_\gamma(s, z, y) q(s, z, y) ds
\le C_1\left(1\wedge \frac{\delta_D(y)}{t^{1/\alpha}}\right)^\gamma
\int^{t/2}_0 \left(1\wedge \frac{\delta_D(z)}{s^{1/\alpha}}\right)^\gamma
q(s, z, y) ds.
\end{equation}
\end{lemma}

\pf
The inequality holds trivially when $\gamma =0$ with $C_1=1$
so for the rest of the proof,  we  assume $\gamma \in (0, 2\alpha)$.
The inequality  \eqref{e:ineqmain1} is obvious if $\delta_D(y) \ge  t^{1/\alpha} $
or $\delta_D(z) \le 2\delta_D(y)$.
So we will assume $\delta_D(y) <   t^{1/\alpha} \wedge (\delta_D(z)/2)$ throughout this proof.
Note that in this case,
\begin{equation}\label{e:tyu1}
|z-y|\geq \delta_D(z)-\delta_D (y) \geq \frac{\delta_D(z)}2\geq
\delta_D(y).
\end{equation}
By \eqref{e:elineq2}, we have
\begin{eqnarray}
&&\int_{(t/2)\wedge |z-y|^\alpha}^{t/2}\left(1\wedge
\frac{\delta_D(z)}{t^{1/\alpha}}\right)^\gamma\psi_\gamma(s, z, y)
q(s, z, y) ds\nonumber\\ &\le& 2^{2\gamma}\left(1\wedge
\frac{\delta_D(y)}{t^{1/\alpha}}\right)^\gamma\int_{(t/2)\wedge
|z-y|^\alpha}^{t/2}\left(1\wedge
\frac{\delta_D(z)}{s^{1/\alpha}}\right)^\gamma q(s, z, y)
ds,     \label{e:e2}
\end{eqnarray}
while by \eqref{e:elineq1}
\begin{align}
& \left(1\wedge
\frac{\delta_D(z)}{t^{1/\alpha}}\right)^\gamma
\int^{(t/2)\wedge |z-y|^\alpha}_0\psi_\gamma(s, z, y)
q(s, z, y)ds            \nonumber      \\
\le&\left( \frac{\delta_D(y)}{t^{1/\alpha}}\right)^\gamma
\left(\frac{\delta_D(z)\wedge
t^{1/\alpha}}{\delta_D(y)}\right)^\gamma \int^{(t/2)\wedge
|z-y|^\alpha}_0     \left(1\wedge
\frac{\delta_D(y)}{s^{1/\alpha}}\right)^\gamma q(s, z, y)ds .
\label{e:new}
\end{align}
In view of
\eqref{e:ineq1}, \eqref{e:tyu1} and \eqref{e:new},
\begin{align}
&\int^{(t/2)\wedge |z-y|^\alpha}_0 \left(1\wedge
\frac{\delta_D(y)}{s^{1/\alpha}}\right)^\gamma
q(s, z, y)ds   \nonumber \\
\asymp &    \int_0^{(t/2)\wedge \delta_D(y)^\alpha}
 \frac{s}{|z-y|^{d+\alpha}} ds
+ \int_{(t/2)\wedge \delta_D(y)^\alpha}^{(t/2)\wedge |z-y|^\alpha}
\left(\frac{\delta_D(y)}{s^{1/\alpha}}\right)^\gamma \frac{s}{|z-y|^{d+\alpha}} ds   \nonumber \\
\asymp& \frac1{|z-y|^{d+\alpha}}
 \left( \left((t/2)\wedge \delta_D(y)^\alpha \right)^2 +
\delta_D(y)^\gamma \left( \left((t/2)\wedge |z-y|^\alpha\right)^{2-\gamma/\alpha} -
\left((t/2)\wedge \delta_D(y)^\alpha \right)^{2-\gamma/\alpha} \right)  \right) \nonumber \\
 \asymp& \frac1{|z-y|^{d+\alpha}}
  \Big( \left((t/2)\wedge \delta_D(y)^\alpha \right)^2 \nonumber \\
  \quad&+
(\delta_D (y)\wedge (t/2)^{1/\alpha})^\gamma
\left( \left((t/2)\wedge |z-y|^\alpha\right)^{2-\gamma/\alpha} -
\left((t/2)\wedge \delta_D(y)^\alpha \right)^{2-\gamma/\alpha} \right)  \Big) \nonumber \\
\asymp& \frac{ (\delta_D (y)\wedge (t/2)^{1/\alpha})^\gamma \,
\left((t/2)\wedge |z-y|^\alpha\right)^{2-\gamma/\alpha}}
{|z-y|^{d+\alpha}}
\,
\asymp\, \frac{ \delta_D (y)^\gamma \,
\left((t/2)\wedge |z-y|^\alpha\right)^{2-\gamma/\alpha}}
{|z-y|^{d+\alpha}} . \label{e:2.8}
\end{align}
On the other hand, using \eqref{e:tyu1} we have
\begin{align}
&\int^{(t/2)\wedge |z-y|^\alpha}_0 \left(1\wedge
\frac{\delta_D(z)}{s^{1/\alpha}}\right)^\gamma
q(s, z, y)ds   \nonumber \\
\asymp &    \int_0^{(t/2)\wedge (\delta_D(z)/2)^\alpha}
 \frac{s}{|z-y|^{d+\alpha}} ds
+ \int_{(t/2)\wedge (\delta_D(z)/2)^\alpha}^{(t/2)\wedge |z-y|^\alpha} \left(\frac{\delta_D(z)}{s^{1/\alpha}}\right)^\gamma \frac{s}{|z-y|^{d+\alpha}} ds   \nonumber \\
\asymp& \frac1{|z-y|^{d+\alpha}}
\left( \left(\frac{t}2\wedge \left(\frac{\delta_D(z)}2\right)^\alpha \right)^2 +
\delta_D(z)^\gamma \left( \left(\frac{t}2\wedge |z-y|^\alpha\right)^{2-\gamma/\alpha} -
\left(\frac{t}2\wedge \left(\frac{\delta_D(z)}2\right)^\alpha \right)^{2-\gamma/\alpha} \right)  \right) \nonumber \\
 \ge& \frac1{|z-y|^{d+\alpha}}
\left(\left(\frac{t}2\wedge \left(\frac{\delta_D(z)}2\right)^\alpha \right)^2 \right.\nonumber \\ \quad&+
\left.\left( \frac{\delta_D(z)}2 \wedge \left(\frac{t}2\right)^{1/\alpha}\right) ^\gamma
\left( \left(\frac{t}2\wedge |z-y|^\alpha\right)^{2-\gamma/\alpha} -
\left( \frac{t}2 \wedge \left(\frac{\delta_D(z)}2\right)^\alpha \right)^{2-\gamma/\alpha} \right)\right) \nonumber \\
\asymp& \frac{ (\delta_D(z) \wedge t^{1/\alpha})^\gamma \,
\left((t/2)\wedge |z-y|^\alpha\right)^{2-\gamma/\alpha}}
{|z-y|^{d+\alpha}} . \label{e:2.9}
\end{align}
One then deduces from \eqref{e:new}--\eqref{e:2.9} and the assumption
  $\delta_D(y) \leq t^{1/\alpha}$ that
\begin{eqnarray*}
&& \left(1\wedge
\frac{\delta_D(z)}{t^{1/\alpha}}\right)^\gamma
\int^{(t/2)\wedge |z-y|^\alpha}_0\psi_\gamma(s, z, y)
q(s, z, y)ds            \nonumber      \\
&\leq & c_1 \left(1\wedge
\frac{\delta_D(y)}{t^{1/\alpha}}\right)^\gamma
(\delta_D(z) \wedge t^{1/\alpha})^\gamma
\frac{\left((t/2)\wedge |z-y|^\alpha\right)^{2-\gamma/\alpha}}
{|z-y|^{d+\alpha}}\\
&\leq & c_2 \int^{(t/2)\wedge |z-y|^\alpha}_0 \left(1\wedge
\frac{\delta_D(z)}{
s^{1/\alpha}}\right)^\gamma
q(s, z, y)ds.
\end{eqnarray*}
This combining with   \eqref{e:e2}
establishes the inequality \eqref{e:ineqmain1}.
\qed

It follows from \eqref{e:ineq*} and \eqref{e:ineq1} that
for every $0<s<t$, and $x, y, z\in \bR^d$,
\begin{eqnarray}
&&\frac{q(s,x,z) q(t-s, z,y)}{q(t,x,y)}\nonumber\\
& \le&
4^{d+\alpha} \frac{s(t-s)}{t}
\left(\frac{ t^{1/\alpha}+  |x-y|}{ (s^{1/\alpha}+
|x-z|)((t-s)^{1/\alpha} + |y-z|)} \right)^{d+\alpha} \nonumber\\
&\le&
4^{d+\alpha} (s \wedge (t-s))
\left(\frac{ (s+(t-s))^{1/\alpha}+  |x-z|+|y-z| }{
(s^{1/\alpha}+  |x-z|)((t-s)^{1/\alpha} + |y-z|)} \right)^{d+\alpha}\nonumber\\
&\le&
2^{(d+\alpha)(3+1/\alpha)} (s \wedge (t-s))
\left(\frac{1}{(s^{1/\alpha}+|x-z|)^{d+\alpha}}+
\frac{1}{((t-s)^{1/\alpha}+|y-z|)^{d+\alpha}} \right)\nonumber\\
&\le& 2^{(d+\alpha)(3+1/\alpha)}
\left(q(s,x,z)+ q(t-s, z,y) \right).\label{e:ppp}
\end{eqnarray}
(See also \cite{BJ}.)

Now we are ready to prove one form of the 3P inequality.
Note that the right hand side of the 3P inequality below has the term
$q(s, x, z)+q(s, z,y)$ rather than
$q(t-s, x, z)+q(s, z,y)$.

\begin{lemma}[3P inequality]\label{l:3p}
For every
$\gamma\in [0, \alpha)$,
there exists a constant $C_2:=
C_2(d, \alpha, \gamma)>0$ such that for all
$(t, x,y, z) \in (0, \infty) \times D\times D \times D$,
\begin{eqnarray*}
\int^t_0\frac{\psi_\gamma(t-s, x, z)q(t-s, x, z)
\psi_\gamma(s, z,y)  q(s, z,y)}{ \psi_\gamma(t, x, y)
q(t, x, y)}ds
\le   C_2\, \int^{t}_0 \left(1\wedge \frac{\delta_D(z)}
{s^{1/\alpha}}\right)^\gamma (q(s, x, z)+q(s, z,y))ds.
\end{eqnarray*}
\end{lemma}

\pf
When $\gamma =0$, the desired inequality follows from
\eqref{e:ppp} with $C_2=2^{(d+\alpha)(3+1/\alpha)}$.
So for the rest of the proof, we assume
$\gamma \in (0, \alpha)$.
 Let
$$
J(t, x, y, z):= \int^t_0\psi_\gamma(t-s, x, z)
q(t-s, x, z) \psi_\gamma(s, z,y) q(s, z,y) ds.
$$
Since
\begin{eqnarray*}
J(t, x, y, z)
 &\le & c_1\left(1\wedge \frac{\delta_D(x)}
{t^{1/\alpha}}\right)^\gamma    \left(1\wedge
\frac{\delta_D(z)}{t^{1/\alpha}}\right)^\gamma  q(t, x, z) \int^{t/2}_0
 \psi_\gamma(s, z, y) q(s, z, y)
ds   \\
&& +  c_1\left(1\wedge \frac{\delta_D(y)}{t^{1/\alpha}}\right)^\gamma
 \left(1\wedge \frac{\delta_D(z)}{t^{1/\alpha}}\right)^\gamma
 q(t, z, y) \int^{t}_{t/2}\psi_\gamma(t-s,
x, z) q(t-s, x,z) ds,
\end{eqnarray*}
we have by Lemma \ref{l:4-1}  that
\begin{eqnarray*}
J(t, x, y, z)
&\le& c_2 \psi_\gamma(t, x, y)
\int^{t/2}_0\left(1\wedge
\frac{\delta_D(z)}{s^{1/\alpha}}\right)^\gamma q(t-s, x, z) q(s, z, y)
ds   \\
&&+ c_2\psi_\gamma(t, x, y)  \int^{t}_{t/2}\left(1\wedge
\frac{\delta_D(z)}{(t-s)^{1/\alpha}}\right)^\gamma q(t-s, x, z) q(s, z, y)
ds  .
\end{eqnarray*}
It then follows from  \eqref{e:ppp} that
\begin{eqnarray*}
J(t, x, y, z)
&\le& c_3 \psi_\gamma(t, x, y) q(t, x, y) \int^{t/2}_0
\left(1\wedge \frac{\delta_D(z)}{s^{1/\alpha}}\right)^\gamma
(q(t-s, x, z)+q(s, z, y)) ds   \\
&& + c_3 \psi_\gamma(t, x, y) q(t, x, y)  \int_{t/2}^t\left(1\wedge
\frac{\delta_D(z)}{(t-s)^{1/\alpha}}\right)^\gamma ( q(t-s, x,z)
+q(s, z, y)) ds   \\
&\le& c_4 \psi_\gamma(t, x, y) q(t, x, y) \int^{t}_0 \left(1\wedge
\frac{\delta_D(z)}{s^{1/\alpha}}\right)^\gamma (q(s, z, y)+q(s, x,z))    ds.
\end{eqnarray*}
Here in the last inequality, we used the fact that
$$
\int_0^{t/2}  \left(1\wedge \frac{\delta_D(z)}{s^{1/\alpha}}\right)^\gamma
q(t-s, x, z) ds \leq c_5 \int_{t/2}^t \left(1\wedge \frac{\delta_D(z)}{s^{1/\alpha}}\right)^\gamma
q(s, x, z)  ds
$$
and
$$
\int_{t/2}^t\left(1\wedge
\frac{\delta_D(z)}{(t-s)^{1/\alpha}}\right)^\gamma  q(s, z, y) ds
\leq c_5 \int_{t/2}^t\left(1\wedge
\frac{\delta_D(z)}{s^{1/\alpha}}\right)^\gamma  q(s, z, y) ds.
$$
The above two inequalities can be easily verified by using the facts
that $q(s, x, y)\asymp q(t, x, y)$ for $s\in [t/2, t]$ and
that
$$
\int_0^{t/2} \left(1\wedge \frac{\delta_D(z)}{s^{1/\alpha}}\right)^\gamma ds
\leq \frac{\alpha}{\alpha-\gamma} 2^{\gamma/\alpha-1}  t \left(1\wedge \frac{\delta_D(z)}{t^{1/\alpha}}\right)^\gamma
\leq c_6 \int_{t/2}^t \left(1\wedge \frac{\delta_D(z)}{s^{1/\alpha}}\right)^\gamma ds,
$$
which follows easily from the assumption
$\gamma \in (0, \alpha)$
by a direct calculation. This completes the proof of the lemma.
\qed

The above 3P inequality immediately implies the following theorem,
which will be used later.

\begin{thm}\label{t:key}
For every
$\gamma\in [0, \alpha)$,
there exists a constant $C_3=
C_3 (d, \alpha, \gamma)>0$ such that for any measure $\mu$ on
$D$
and any $(t, x,y)\in (0, \infty)\times D \times D$,
\begin{eqnarray*}
&&\int^t_0\int_D\psi_\gamma(t-s, x, z) \psi_\gamma(s, z,y) q(t-s, x,
z) q(s, z,y) \mu(dz)ds\nonumber\\
&\le&C_3\, \psi_\gamma(t, x, y) q(t,x,y) \sup_{u\in
D}\int^t_0\int_{D}\left(1\wedge
\frac{\delta_D(z)}{s^{1/\alpha}}\right)^\gamma q(s,u,z) \mu(dz)ds.
\end{eqnarray*}
\end{thm}

The results of the remainder of this section are geared towards
dealing with the discontinuous part of $A^{\mu, F}$.

\begin{lemma}\label{l:4-2}
For every
$\gamma \in [0, 2\alpha)$,
there exists a constant
$C_4:=C_4(d, \alpha, \gamma)>1$ such that for all $(t, y,z,w)\in (0,
\infty)\times D \times D \times D$,
\begin{align}
& \left(1\wedge
\frac{\delta_D(z)}{t^{1/\alpha}}\right)^\gamma\int^{t/2}_0
\psi_\gamma(s, w, y)
 q(s, w, y) ds\nonumber \\
\le& C_4\left(1\wedge
\frac{\delta_D(y)}{t^{1/\alpha}}\right)^\gamma \left(1
  + \frac{|y-z| \wedge |z-w| \wedge t^{1/\alpha}}{|y-w|}
\right)^\gamma\int^{t/2}_0 \left(1\wedge
\frac{\delta_D(w)}{s^{1/\alpha}}\right)^\gamma
q(s, w, y) ds.\label{e:ineqmain2}
\end{align}
\end{lemma}

\pf
The desired inequality holds trivially for $\gamma =0$ with $C_4=1$
so for the rest of the proof we assume $\gamma \in (0, 2\alpha)$.
The inequality \eqref{e:ineqmain2} is obvious if $\delta_D(y) \ge  t^{1/\alpha}$ or $\delta_D(z) \le 2\delta_D(y)$,
so we will assume $\delta_D(y) <  t^{1/\alpha} \wedge (\delta_D(z)/2)$ in the remainder of this proof. Note that in this case
\begin{equation}\label{e:2.14}
|y-z|\geq \delta_D(z)-\delta_D(y) \geq \frac{\delta_D(z)}2 \geq \delta_D (y),
\end{equation}
By \eqref{e:elineq1}, \eqref{e:ineq*} and our assumption $\delta_D(y) <  t^{1/\alpha}$, we have that
\begin{equation}\label{e:2.13}
\left( 1 \wedge \frac{\delta_D(z)}{t^{1/\alpha}}\right)
\left( 1 \wedge \frac{\delta_D(y)}{s^{1/\alpha}}\right)
\leq 2  \, \frac{\delta_D(y)}{t^{1/\alpha}}
\left( \frac{\delta_D (z)\wedge t^{1/\alpha}}
{s^{1/\alpha}+\delta_D (y)} \right) =2 \left( 1 \wedge \frac{\delta_D(y)}{t^{1/\alpha}}\right)
\left( \frac{\delta_D (z)\wedge t^{1/\alpha}}
{s^{1/\alpha}+\delta_D (y)}\right).
\end{equation}
When $s\geq |y-w|^\alpha$, by \eqref{e:2.14},
\begin{equation}
  \frac{\delta_D (z)\wedge t^{1/\alpha}}{s^{1/\alpha}+\delta_D (y)}
 \leq   2\, \frac{|y-z|\wedge t^{1/\alpha}}{|y-w|}
 \leq 2 \left( 1+ \frac{|y-z|\wedge |z-w|
 \wedge t^{1/\alpha}}{|y-w|}
 \right), \nonumber
\end{equation}
where the last inequality is due to the
fact
 $|y-z| \le |y-w| + ( |y-z| \wedge |z-w|)$.
This together with
\eqref{e:2.13} implies that
\begin{align}
& \left(1\wedge
\frac{\delta_D(z)}{t^{1/\alpha}}\right)^\gamma
\int^{ t/2 }_{(t/2)\wedge |y-w|^\alpha}
\left(1\wedge \frac{\delta_D(y)}{s^{1/\alpha}}\right)^\gamma
\left(1\wedge \frac{\delta_D(w)}{s^{1/\alpha}}\right)^\gamma
 q(s, w, y) ds\nonumber \\
\le& 4^\gamma \left(1 \wedge \frac{\delta_D(y)}{t^{1/\alpha}}\right)^\gamma
\left(1 + \frac{|y-z| \wedge |z-w| \wedge t^{1/\alpha}}{|y-w|}
\right)^\gamma
\int^{ t/2 }_{(t/2)\wedge |y-w|^\alpha}\left(1\wedge \frac{\delta_D(w)}{s^{1/\alpha}}\right)^\gamma q(s, w, y) ds.\label{e:2.16}
\end{align}
 On the other hand, by \eqref{e:2.14},
\begin{eqnarray}
&&  \int_0^{(t/2)\wedge |y-w|^\alpha}
\left( \frac{\delta_D (z)\wedge t^{1/\alpha}}
{s^{1/\alpha}+\delta_D (y)}\right)^\gamma \left(1\wedge \frac{\delta_D(w)}{s^{1/\alpha}}\right)^\gamma q(s, w, y) ds  \nn \\
&\leq &   2^\gamma  \int_0^{(t/2)\wedge |y-w|^\alpha}
\left( \frac{|y-z|\wedge t^{1/\alpha}}
{s^{1/\alpha} }\right)^\gamma \left(1\wedge \frac{\delta_D(w)}{s^{1/\alpha}}\right)^\gamma\, \frac{s}{|y-w|^{d+\alpha}} ds \nn \\
&= & c_1 \, \frac{\left( |y-z|\wedge t^{1/\alpha}\right)^\gamma}{|y-w|^{d+\alpha}}
\int_0^{(t/2)\wedge |y-w|^\alpha}
{s^{1-\gamma/\alpha} } \left(1\wedge \frac{\delta_D(w)}{s^{1/\alpha}}\right)^\gamma ds.
\label{e:2.15}
\end{eqnarray}
We claim that
\begin{align}
 \int_0^{(t/2)\wedge |y-w|^\alpha}
{s^{1-\gamma/\alpha} } \left(1\wedge \frac{\delta_D(w)}{s^{1/\alpha}}\right)^\gamma ds
\asymp   \left(\frac{t}2\wedge |y-w|^\alpha\right)^{-\gamma/\alpha}
\int_0^{(t/2)\wedge |y-w|^\alpha}
{s } \left(1\wedge \frac{\delta_D(w)}{s^{1/\alpha}}\right)^\gamma ds . \label{e:2.811}
 \end{align}
The case $\delta_D(w) > (t/2)^{1/\alpha}$ is clear. If $\delta_D(w) \le  |y-w|\wedge (t/2)^{1/\alpha} $,
\begin{align}
&\int_0^{(t/2)\wedge |y-w|^\alpha}
{s^{1-\gamma/\alpha} } \left(1\wedge \frac{\delta_D(w)}{s^{1/\alpha}}\right)^\gamma ds   \nonumber \\
= &    \int_0^{\delta_D(w)^\alpha}
 s^{1-\gamma/\alpha} ds
+ \delta_D(w)^\gamma\int_{\delta_D(w)^\alpha}^{(t/2)\wedge |y-w|^\alpha}
 s^{1-2\gamma/\alpha} ds   \nonumber \\
\asymp& \
  \delta_D(w)^{2\alpha-\gamma} +
\delta_D(w)^\gamma
\left( \left(\frac{t}2\wedge |y-w|^\alpha\right)^{2-2\gamma/\alpha} - \delta_D(w)^{2(\alpha-\gamma)} \right)   \nonumber \\
\asymp&\delta_D(w)^\gamma
\left(\frac{t}2\wedge |y-w|^\alpha\right)^{2-2\gamma/\alpha}\nonumber \\
=&\left(\frac{t}2\wedge |y-w|^\alpha\right)
^{-\gamma/\alpha} \left(\delta_D(w)\right)^\gamma
\left(\frac{t}2\wedge |y-w|^\alpha\right)^{2-\gamma/\alpha} \nonumber \\
\asymp&\left(\frac{t}2\wedge |y-w|^\alpha\right)
^{-\gamma/\alpha} \int_0^{(t/2)\wedge |y-w|^\alpha}
{s } \left(1\wedge \frac{\delta_D(w)}{s^{1/\alpha}}\right)^\gamma ds
 . \nonumber
\end{align}
The remaining case $ |y-w| < \delta_D(w)\le (t/2)^{1/\alpha}$ is simpler.
Thus we have proved the claim \eqref{e:2.811}. Now by \eqref{e:2.15} and
\eqref{e:2.811},
 \begin{eqnarray*}
&&  \int_0^{(t/2)\wedge |y-w|^\alpha}
\left( \frac{\delta_D (z)\wedge t^{1/\alpha}}
{s^{1/\alpha}+\delta_D (y)}\right)^\gamma \left(1\wedge \frac{\delta_D(w)}{s^{1/\alpha}}\right)^\gamma q(s, w, y) ds \\
&\leq & c_2 \left( \frac{ |y-z|\wedge t^{1/\alpha}}{|y-w|\wedge t^{1/\alpha}}\right)^\gamma  \int_0^{(t/2)\wedge |y-w|^\alpha}
\left(1\wedge \frac{\delta_D(w)}{s^{1/\alpha}}\right)^\gamma  \frac{s}{|y-w|^{d+\alpha}} ds \\
&\leq & c_2 \left( 1+ \frac{ |y-z|\wedge t^{1/\alpha}}{|y-w| }\right)^\gamma  \int_0^{(t/2)\wedge |y-w|^\alpha}
 \left(1\wedge \frac{\delta_D(w)}{s^{1/\alpha}}\right)^\gamma q(s, w, y) ds \\
&\leq & 2c_2 \left( 1+ \frac{ |y-z|\wedge |z-w|\wedge t^{1/\alpha}}
{|y-w| }\right)^\gamma  \int_0^{(t/2)\wedge |y-w|^\alpha}
 \left(1\wedge \frac{\delta_D(w)}{s^{1/\alpha}}\right)^\gamma q(s, w, y) ds .
\end{eqnarray*}
Here again the last inequality is due to the fact that
 $|y-z| \le |y-w| + ( |y-z| \wedge |z-w|)$.
This together with \eqref{e:2.13} and \eqref{e:2.16} establishes
the inequality \eqref{e:ineqmain2}.
\qed
\

In the remainder of this section.
we use the following notation: For any $(x, y)\in D\times D$,
\begin{eqnarray*}
V_{x,y}&:=&{\{(z,w)\in D\times D:
|x-y| \ge 4 (|y-w| \wedge |x-z|)
  \}},\\
U_{x,y}&:=&(D \times D) \setminus
V_{x,y}.
\end{eqnarray*}
Recall that, for any bounded function $F$ on $D\times D$ we use $\|F\|_\infty$
to denote $\|F\|_{L^\infty (D \times D)}$.

Now we are ready to prove  the following generalized 3P inequality.

\begin{thm}[Generalized 3P inequality]\label{t:G3p}
For every
$\gamma \in [0, \alpha \wedge d)$,
there exists a constant
$C_5:=C_5(\alpha, \gamma, d)>0$ such that
for any nonnegative  bounded function $F(x,y)$ on $D\times D$, the following are true for
$(t, x, y)\in (0, \infty) \times D\times D$.

\noindent (a)
If $|x-y| \le t^{1/\alpha}$, then
\begin{eqnarray*}
&&\int^{t}_{0}\int_{D\times D}   \frac{\psi_\gamma(t-s, x, z)q(t-s,
x, z)
\psi_\gamma(s, w,y)  q(s, w,y)}
{ \psi_\gamma(t, x, y) q(t, x,
y)} \frac{F(z,w)}{|z-w|^{d+\alpha}} dzdwds \nonumber\\
&\le& C_5\int^{t}_{0}\int_{D\times D}   \left(1\wedge \frac{\delta_D(z)}{s^{1/\alpha}}\right)^\gamma
{q(s, x,z)}\left(1 + \frac{ |z-w|
\wedge t^{1/\alpha}}{|x-z|} \right)^\gamma
\frac{F(z,w)}{|z-w|^{d+\alpha}} dzdwds \nonumber\\
&&+ C_5\int^{t}_{0}\int_{D\times D}   \left(1\wedge \frac{\delta_D(w)}{s^{1/\alpha}}\right)^\gamma
{q(s, y,w)}\left(1 + \frac{ |z-w|
\wedge t^{1/\alpha}}{|y-w|} \right)^\gamma
\frac{F(z,w)}{|z-w|^{d+\alpha}} dzdwds.
\end{eqnarray*}

\noindent (b)
If $|x-y| > t^{1/\alpha}$, then
\begin{eqnarray*}
&&\int^{t}_{0}\int_{U_{x,y}}   \frac{\psi_\gamma(t-s, x, z)q(t-s,
x, z)\psi_\gamma(s, w,y)  q(s, w,y)}
{ \psi_\gamma(t, x, y) q(t, x,
y)} \frac{F(z,w)}{|z-w|^{d+\alpha}} dzdwds \nonumber\\
&\le& C_5\int^{t}_{0}\int_{U_{x,y}}   \left(1\wedge \frac{\delta_D(z)}{s^{1/\alpha}}\right)^\gamma
{q(s, x,z)}\left(1 + \frac{ |z-w|
\wedge t^{1/\alpha}}{|x-z|} \right)^\gamma
\frac{F(z,w)}{|z-w|^{d+\alpha}} dzdwds \nonumber\\
&&+ C_5\int^{t}_{0}\int_{U_{x,y}}  \left(1\wedge \frac{\delta_D(w)}{s^{1/\alpha}}\right)^\gamma
{q(s, y,w)}\left(1 + \frac{ |z-w|
\wedge t^{1/\alpha}}{|y-w|} \right)^\gamma
\frac{F(z,w)}{|z-w|^{d+\alpha}} dzdwds .
\end{eqnarray*}

\noindent (c)
If $|x-y| > t^{1/\alpha}$, then
$$
\int^{t}_{0}\int_{V_{x,y}}   \frac{\psi_\gamma(t-s, x, z)q(t-s,
x, z) \psi_\gamma(s, w,y)  q(s, w,y)}
{ \psi_\gamma(t, x, y) q(t, x,
y)} \frac{F(z,w)}{|z-w|^{d+\alpha}} dzdwds \, \le \,  C_5 \|F\|_{\infty }.
$$
\end{thm}

\pf
By Lemma \ref{l:4-2}, we get that
\begin{align}
& \int^t_0\int_{D\times D}
\frac{
\psi_\gamma(t-s, x, z)
q(t-s, x, z) \psi_\gamma(s, w,y) q(s, w,y)}
{\psi_\gamma(t, x, y)}
\frac{F(z,w)}{|z-w|^{d+\alpha}} dzdwds   \nonumber\\
\le& c_1 \int_{D\times D} \int^{t/2}_{0} \left(1\wedge \frac{\delta_D(w)}{s^{1/\alpha}}\right)^\gamma q(s,
w,y)q(t-s, x,z)\left(1 + \frac{ |z-w| \wedge t^{1/\alpha}}{|y-w|}
\right)^\gamma
ds \frac{F(z,w)}{|z-w|^{d+\alpha}} dzdw\nonumber\\
&+ c_1\int_{D\times D} \int^{t}_{t/2} \left(1\wedge \frac{\delta_D(z)}{(t-s)^{1/\alpha}}\right)^\gamma q(s, w,y)
q(t-s, x,z)\left(1 + \frac{|z-w|\wedge t^{1/\alpha}}{|x-z|}
\right)^\gamma
ds \frac{F(z,w)}{|z-w|^{d+\alpha}} dzdw. \label{e:newe1}
\end{align}
If $|x-y| \le t^{1/\alpha}$ and $s\in (0, t/2]$, we have $ q(t-s, x,z)
\le 2^{d/\alpha}q(t,x,y)$, and if $|x-y| \le t^{1/\alpha}$ and $s\in (t/2, t]$, we have $ q(s, w,y)
\le 2^{d/\alpha}q(t,x,y)$.
Thus (a) follows immediately from \eqref{e:newe1}.

In the remainder of this proof,
we fix $(t,x,y)\in (0, \infty) \times D \times D$ with $|x-y|> t^{1/\alpha}$.
Let
\begin{eqnarray*}
U_1&:=&{\{(z,w)\in D\times D: |y-w| >
4^{-1}|x-y|,
|y-w| \ge |x-z|  \}}, \\
U_2&:=&{\{(z,w)\in D\times D: |x-z| >
4^{-1}|x-y| \}}.
\end{eqnarray*}
Since
$ q(t-s, x,z) \le 4^{d+\alpha} q(t,x,y)$
for $(s, z, w)\in (0, t)\times U_2$,
by Lemma \ref{l:4-2}, we have
\begin{align}
& \int^{t/2}_0\int_{U_2}
\frac{\psi_\gamma(t-s, x, z)q(t-s, x, z) \psi_\gamma(s, w,y) q(s, w,y)}{\psi_\gamma(t, x, y)}
\frac{F(z,w)}{|z-w|^{d+\alpha}} dzdwds \nonumber\\
\le& c_2  \int_{U_2} \int^{t/2}_{0} \left(1\wedge \frac{\delta_D(w)}{s^{1/\alpha}}\right)^\gamma q(s,
w,y)q(t-s, x,z)\left(1 + \frac{ |z-w| \wedge t^{1/\alpha}}{|y-w|}
\right)^\gamma
ds \frac{F(z,w)}{|z-w|^{d+\alpha}} dzdw\nonumber\\
\le& c_3  q(t, x,y)\int_{U_2} \int^{t/2}_{0} \left(1\wedge \frac{\delta_D(w)}{s^{1/\alpha}}\right)^\gamma q(s,
w,y)\left(1 + \frac{ |z-w| \wedge t^{1/\alpha}}{|y-w|}
\right)^\gamma
ds \frac{F(z,w)}{|z-w|^{d+\alpha}} dzdw\label{e:new_1}
\end{align}
and, similarly
\begin{align}
& \int^t_{t/2}\int_{U_1}
\frac{\psi_\gamma(t-s, x, z)
q(t-s, x, z) \psi_\gamma(s, w,y) q(s, w,y)}{\psi_\gamma(t, x, y)}
\frac{F(z,w)}{|z-w|^{d+\alpha}} dzdwds\nonumber\\
\le&  c_4\int_{U_1} \int^{t}_{t/2} \left(1\wedge \frac{\delta_D(z)}{(t-s)^{1/\alpha}}\right)^\gamma q(s, w,y)
q(t-s, x,z)\left(1 + \frac{|z-w|\wedge t^{1/\alpha}}{|x-z|}
\right)^\gamma
ds \frac{F(z,w)}{|z-w|^{d+\alpha}} dzdw\nonumber\\
\le&  c_5  q(t, x,y) \int_{U_1} \int^{t}_{t/2} \left(1\wedge \frac{\delta_D(z)}{(t-s)^{1/\alpha}}\right)^\gamma
q(t-s, x,z)\left(1 + \frac{|z-w|\wedge t^{1/\alpha}}{|x-z|}
\right)^\gamma
ds \frac{F(z,w)}{|z-w|^{d+\alpha}} dzdw. \label{e:new_2}
\end{align}

On the other hand, we observe that, since  $ q(s, w,y) \le 4^{d+\alpha} q(t,x,y)$
for $(s, z, w)\in (0, t/2]\times U_1$,
\begin{align*}
& \int^{t/2}_0\int_{U_1}\psi_\gamma(t-s, x, z)
q(t-s, x, z) \psi_\gamma(s, w,y) q(s, w,y)
\frac{F(z,w)}{|z-w|^{d+\alpha}} dzdwds\\
\le&
c_6  \psi_\gamma(t, x, z)q(t, x, y)
\int_{U_1}q(t, x, z)
\int^{t/2}_0 \psi_\gamma(s, w,y) ds
\frac{F(z,w)}{|z-w|^{d+\alpha}} dzdw.
\end{align*}
Now, applying the inequality
$$
\int^{t/2}_0 \psi_\gamma(s, w,y) ds \le
\int^{t/2}_0 \left(1\wedge \frac{\delta_D(y)}{s^{1/\alpha}}\right)^\gamma ds
\le
\frac{\alpha}{\alpha-\gamma} 2^{\gamma/\alpha-1}
 t\left(1\wedge \frac{\delta_D(y)}{t^{1/\alpha}}\right)^\gamma,
$$
we get
\begin{align}
& \int^{t/2}_0\int_{U_1}\psi_\gamma(t-s, x, z)
q(t-s, x, z) \psi_\gamma(s, w,y) q(s, w,y)
\frac{F(z,w)}{|z-w|^{d+\alpha}} dzdwds\nonumber\\
\le&
c_7  \psi_\gamma(t, x, y)q(t, x, y)
\int_{U_1}q(t, x, z)
t\left(1\wedge \frac{\delta_D(z)}{t^{1/\alpha}}\right)^\gamma
\frac{F(z,w)}{|z-w|^{d+\alpha}} dzdw\nonumber\\
\le&
c_8  \psi_\gamma(t, x, y)q(t, x, y)
\int_{U_1}
\int_0^{t/2} q(t-s, x, z)\left(1\wedge \frac{\delta_D(z)}{(t-s)^{1/\alpha}}\right)^\gamma
\frac{F(z,w)}{|z-w|^{d+\alpha}}ds dzdw.\label{e:new_3}
\end{align}
Similarly
\begin{align}
& \int^t_{t/2}\int_{U_2}\psi_\gamma(t-s, x, z)
q(t-s, x, z) \psi_\gamma(s, w,y) q(s, w,y)
\frac{F(z,w)}{|z-w|^{d+\alpha}} dzdwds\nonumber\\
\le&  c_{9}\psi_\gamma(t, x, y)q(t, x, y)\int_{U_2} \int^{t}_{t/2}
\left(1\wedge \frac{\delta_D(w)}{s^{1/\alpha}}\right)^\gamma q(s, w,y)
\frac{F(z,w)}{|z-w|^{d+\alpha}}ds dzdw. \label{e:new_4}
\end{align}
Since $U_{x, y}=U_1\cup U_2$, from \eqref{e:new_1}--\eqref{e:new_4}, we know that (b) is true.

Note that for $(z,w)\in V_{x,y}$, we have
$|z-w| \ge |x-y| -(|x-z|+|y-w|) \ge 2^{-1} |x-y|.$
Thus, by Lemma \ref{l:4-2} and \eqref{e:hkein1}, it is easy to see that
\begin{align*}
& \int^t_0\int_{V_{x,y}}\frac{\psi_\gamma(t-s, x, z)
q(t-s, x, z) \psi_\gamma(s, w,y) q(s, w,y)}{\psi_\gamma(t, x, y)}
\frac{F(z,w)}{q(t,x,y)|z-w|^{d+\alpha}} dzdwds\\
\le& c_{10} \|F\|_\infty   t^{-1}\int_{V_{x,y}} \int^{t/2}_{0} {q(s,
w,y)q(t-s, x,z)}\left(1 +\frac{ t^{1/\alpha}}{|y-w|}
\right)^\gamma dsdzdw\\
&+ c_{10} \|F\|_\infty  t^{-1} \int_{V_{x,y}} \int^{t}_{t/2} {q(s,
w,y)q(t-s, x,z)}\left(1 +\frac{ t^{1/\alpha}}{|x-z|}
\right)^\gamma dsdzdw\\
\le &c_{11} \|F\|_\infty   t^{-1} \int^{t}_{0}\left(\int_{D} {q(s,
w,y)}\left(1 +\frac{ t^{1/\alpha}}{|y-w|}
\right)^\gamma dw  +   \int_{D} {q(s,
x,z)}\left(1 +\frac{ t^{1/\alpha}}{|x-z|}
\right)^\gamma  dz \right) ds.
\end{align*}
Since, using $\gamma \in (0, \alpha \wedge d)$,
\begin{eqnarray*}
&&\int^{t}_{0}\left(\int_{D} {q(s,
w,y)}\left(1 +\frac{ t^{1/\alpha}}{|y-w|}
\right)^\gamma dw  +   \int_{D} {q(s,
x,z)}\left(1 +\frac{ t^{1/\alpha}}{|x-z|}
\right)^\gamma  dz \right) ds\\
&\le & 2^{d+\alpha+1}\int_0^t
\int_{\bR^d}   \frac{s}{(s^{1/\alpha}+|w|)^{d+\alpha}}
\left(\frac{ t^{1/\alpha}}{|w|} \right)^\gamma dwds\\
&= & c_{12}\left(\int^\infty_{0}  \frac{u^{d-1-\gamma}
du}{(1+u)^{d+\alpha}}  \right)  t^{\gamma /\alpha} \int_0^t
s^{-\gamma /\alpha} ds \le c_{13} t,
\end{eqnarray*}
(c) follows immediately.
\qed

\section{Heat kernel estimates}\label{s:3}

In this section we give the proof of our main result, Theorem \ref{t:main}.
Throughout this section,
we fix $\gamma\in [0, \alpha\wedge d)$.
Recall the definition of $p^k(t, x, y)$ given by \eqref{e:i0}.

Using \eqref{e:hke}, \eqref{e:1.5a},
Theorems \ref{t:key} and \ref{t:G3p},
we can choose a  constant
\begin{equation}\label{e:M1}
M= M (\alpha, \gamma, d, C_0)>
\frac{\alpha}{\alpha-\gamma} 2^{2\gamma/\alpha+d+\alpha+1} C_0^4
(C_1 \vee C_4)
\end{equation}
such that for any
 $\mu$ in ${\bf K}_{\alpha, \gamma}$, any
 measurable function  $F$ with $F_1=e^F-1 \in {\bf J}_{\alpha, \gamma}$
and any $(t, x,y)\in
(0, 1]\times D
\times D$,
\begin{eqnarray}
&&\int^t_0\int_D\psi_\gamma(t-s, x, z) \psi_\gamma(s, z,y) q(t-s, x,
z) q(s, z,y) |\mu|(dz)ds\nonumber\\
&\le&M\, p_D(t, x, y) N_\mu^{\alpha, \gamma}(t),\label{e:M2}\\
&& \int^{t}_{0}\int_{D\times D}   {\psi_\gamma(t-s, x,
z)q(t-s, x, z)
\psi_\gamma(s,w,y)
q(s,w,y)}\frac{c(z, w) |F_1|(z,w)}{|z-w|^{d+\alpha}} dzdwds \nonumber \\
&\le& M\, p_D(t, x, y)(N_{F_1}^{\alpha, \gamma}(t)+
\|F_1\|_{\infty}\text{{\bf 1}}_{\{    |x-y| > t^{1/\alpha}
\}})\label{e:M3}
\end{eqnarray}
and
\begin{eqnarray}
&& \int^{t}_{0}\int_{
U_{x,y}}  {\psi_\gamma(t-s, x, z)q(t-s, x, z)
\psi_\gamma(s,w,y)
q(s,w,y)}
\frac{c(z, w)|F_1|(z,w)}{|z-w|^{d+\alpha}} dzdwds \nonumber\\
&\le& M\, p_D(t, x, y)N_{F_1}^{\alpha, \gamma}(t). \label{e:M4}
\end{eqnarray}

In the remainder of this section, we fix
a signed measure $\mu\in{\bf K}_{\alpha, \gamma}$,
a measurable function  $F$ with $F_1=e^F-1 \in {\bf J}_{\alpha, \gamma}$ and the constant $M>0$ in \eqref{e:M1}.

\begin{lemma}\label{l:i1}
For every $k \ge 0$
and $(t, x)\in (0, 1]\times D$,
\begin{eqnarray}\label{e:i1}
\int_D|p^k(t, x, y)|dy \le C_0^2M^k
\,\left(1\wedge
\frac{\delta_D(x)}{t^{1/\alpha}}\right)^\gamma \left(N^{\alpha,
\gamma}_{\mu, F_1}(t)\right)^k.
\end{eqnarray}
\end{lemma}

\pf We use induction on
$k \ge 0$. By \eqref{e:hkein2}, \eqref{e:i1} is clear when $k=0$.
Suppose \eqref{e:i1} is true for $k-1\ge 0$.
Then
by \eqref{e:i0} we have
\begin{align*}
\int_D p^k(t, x, y) dy
=&\int^{t/2}_0 \left(\int_{D}p^0(t-s, x, z)\Big(\int_D
p^{k-1}(s, z, y)dy\Big)\mu(dz) \right)ds \\
&+\int^{t/2}_0 \left(\int_{D}\int_{D} p^0(t-s, x,
z)\frac{c(z, w)F_1(z,w)}{|z-w|^{d+\alpha}} \int_D
p^{k-1}(s, w, y)dy) dzdw \right)ds\\
&+\int^t_{t/2} \left(\int_{D}p^0(t-s, x,
z)\Big(\int_Dp^{k-1}(s, z, y)dy\Big)\mu(dz) \right)ds \\
&+ \int^t_{t/2} \left(\int_{D}\int_{D} p^0(t-s, x,
z)\frac{c(z, w)F_1( z,w)}{|z-w|^{d+\alpha}} \int_D p^{k-1}(s, w, y)dy)
dzdw \right)ds.
\end{align*}
Thus using
\eqref{e:hke} and our induction hypothesis, we have
\begin{align*}
&\int_D |p^k(t, x, y)| dy\\
\le& 2^{\gamma/\alpha} C_0^3 M^{k-1}\big(N^{\alpha, \gamma}_{\mu,
F_1}(t)\big)^{k-1} \left( \left(1\wedge \frac{
\delta_D(x)}{t^{1/\alpha}}\right)^\gamma\int^{t/2}_0 \left(\int_{D}
\left(1\wedge \frac{\delta_D(z)}{(t-s)^{1/\alpha}}\right)^\gamma
q(t-s, x, z)|\mu|(dz) \right)ds \right. \\
& + C_0 \left(1\wedge \frac{
\delta_D(x)}{t^{1/\alpha}}\right)^\gamma\int^{t/2}_0
\left(\int_{D}\int_{D}
\left(1\wedge \frac{\delta_D(z)}{(t-s)^{1/\alpha}}\right)^\gamma
q(t-s, x, z)
\frac{|F_1|(z,w)}
{|z-w|^{d+\alpha}} dzdw \right)ds\\
& + \int_{D}\left(1\wedge
\frac{\delta_D(z)}{t^{1/\alpha}}\right)^\gamma \int^t_{t/2}
\psi_\gamma(t-s, x,
z) q(t-s, x, z)ds |\mu|(dz) \\
& \left. + C_0 \int_{D}\int_{D} \left(1\wedge \frac{\delta_D(w)}{t^{1/\alpha}}
\right)^\gamma
\int^t_{t/2} \psi_\gamma(t-s, x, z)
q(t-s, x,z)
ds \frac{|F_1|(z,w)}{|z-w|^{d+\alpha}} dzdw \right) .
\end{align*}
Applying \eqref{e:M1}, Lemmas \ref{l:4-1} and \ref{l:4-2}, the above is less than
\begin{align*}
& 4^{-1}C_0 M^{k}\big(N^{\alpha, \gamma}_{\mu,
F_1}(t)\big)^{k-1} \left( \left(1\wedge
\frac{\delta_D(x)}{t^{1/\alpha}}\right)^\gamma\int^{t/2}_0
\left(\int_{D} \left(1\wedge
\frac{\delta_D(z)}{(t-s)^{1/\alpha}}\right)^\gamma
q(t-s, x, z)|\mu|(dz) \right)ds  \right. \\
&+ C_0 \left(1\wedge
\frac{\delta_D(x)}{t^{1/\alpha}}\right)^\gamma \int^{t/2}_0
\left(\int_{D}\int_{D}
\left(1\wedge \frac{\delta_D(z)}{(t-s)^{1/\alpha}}\right)^\gamma
q(t-s, x, z)
\frac{|F_1|(z, w)}{|z-w|^{d+\alpha}} dzdw \right)ds\\
&+ \int_{D}
\left(1\wedge \frac{\delta_D(x)}{t^{1/\alpha}} \right)^\gamma
\int^t_{t/2} \left(1\wedge
\frac{\delta_D(z)}{(t-s)^{1/\alpha}}\right)^\gamma
q(t-s, x, z)ds |\mu|(dz) \\
&  + C_0  \int_{D}\int_{D} \left(1\wedge
\frac{\delta_D(x)}{t^{1/\alpha}}\right)^\gamma\left(1 + \frac{|x-w|
\wedge |z-w|\wedge t^{1/\alpha}}{|x-z|}
\right)^\gamma
 \\
 &\left.\quad \times \int^{t}_{t/2}
\left(1\wedge \frac{\delta_D(z)}{(t-s)^{1/\alpha}}\right)^\gamma
q(t-s, x, z)
\frac{|F_1|(z,w)}{|z-w|^{d+\alpha}} dzdw  \right) \\
\le&C_0^2 M^{k}\,\left(1\wedge
\frac{\delta_D(x)}{t^{1/\alpha}}\right)^\gamma  (N^{\alpha,
\gamma}_{\mu, F_1}(t))^k.
\end{align*}
\qed

\begin{lemma}\label{l:i2}
For every $k \ge 0$
and $(t, x, y)\in (0, 1]\times D\times D$,
\begin{eqnarray*}
\int_0^t \int_D p_D(t-s, x, z) dz \int_D
|p^k(s, w, y)| dw ds
\,\le\,t\,
\frac{\alpha}{\alpha-\gamma}2^{2\gamma/\alpha}
\, C_0^4\, M^k\,
 \psi_\gamma (t,x,y)\big(N^{\alpha,
\gamma}_{\mu, F_1}(t)\big)^k.
\end{eqnarray*}
\end{lemma}

\pf By \eqref{e:hke} and Lemma \ref{l:i1},
\begin{eqnarray*}
&&\int_0^t \int_D p_D(t-s, x, z) dz \int_D
|p^k(s, w, y)|
 dw
ds \\
&=& \int_0^{t/2} \int_D p_D(t-s, x, z) dz \int_D
|p^k(s, w, y)| dw
ds+\int_{t/2}^t \int_D p_D(t-s, x, z) dz \int_D
|p^k(s, w, y)| dw
ds\\
&\le&
C_0 \int_0^{t/2} \int_D \left(1\wedge
\frac{\delta_D(x)}{(t-s)^{1/\alpha}}\right)^\gamma   q(t-s, x, z) dz
\int_D
|p^k(s, w, y)| dw ds\\
&&+
C_0 \int_{t/2}^t \int_D\left(1\wedge \frac{\delta_D(x)}
{(t-s)^{1/\alpha}}\right)^\gamma   q(t-s, x, z) dz \int_D
|p^k(s, w, y)| dw ds\\
&\le&
C_0^3 M^{k} \int_0^{t/2} \int_D \left(1\wedge
\frac{\delta_D(x)}{(t-s)^{1/\alpha}}\right)^\gamma   q(t-s, x, z) dz
\,\left(1\wedge \frac{\delta_D(y)}{s^{1/\alpha}}\right)^\gamma
\big(N_{\mu, F_1}(s)\big)^k ds\\
&&+
C_0^3 M^{k} \int_{t/2}^t \int_D\left(1\wedge
\frac{\delta_D(x)} {(t-s)^{1/\alpha}}\right)^\gamma   q(t-s, x, z)
dz \,\left(1\wedge \frac{\delta_D(y)}{s^{1/\alpha}}\right)^\gamma
\big(N_{\mu, F_1}(s)\big)^k ds\\
&\le&2^{\gamma/\alpha}
C_0^3M^{k} \big(N^{\alpha, \gamma}_{\mu,
F_1}(t)\big)^k \left(1\wedge
\frac{\delta_D(x)}{t^{1/\alpha}}\right)^\gamma
\left(\int_0^{t/2}\left(1\wedge
\frac{\delta_D(y)}{s^{1/\alpha}}\right)^\gamma ds \right) \int_D
q(t-s, x, z) dz \\
&&+2^{\gamma/\alpha}
C_0^3M^{k}  \big(N^{\alpha, \gamma}_{\mu,
F_1}(t)\big)^k \left(1\wedge
\frac{\delta_D(y)}{t^{1/\alpha}}\right)^\gamma \int_D
\left(\int_{t/2}^t\left(1\wedge
\frac{\delta_D(x)}{(t-s)^{1/\alpha}}\right)^\gamma  ds \right)q(t-s,
x, z) dz.
\end{eqnarray*}
Using
$$
\int_0^{t/2} \left(1\wedge
\frac{\delta_D(x)}{s^{1/\alpha}}\right)^\gamma ds
\le\frac{\alpha}{\alpha-\gamma}2^{\gamma/\alpha-1}
t \left(1\wedge
\frac{\delta_D(x)}{t^{1/\alpha}}\right)^\gamma,
$$
we get that
\begin{eqnarray*}
&&\int_0^t \int_D p_D(t-s, x, z) dz \int_D p^k(s, w, y) dw
ds \\
&&\,\le\,
\frac{\alpha}{\alpha-\gamma}2^{2\gamma/\alpha}
C_0^3 M^{k} \big(N^{\alpha, \gamma}_{\mu,
F_1}(t)\big)^k \left(1\wedge \frac{\delta_D(y)}
{t^{1/\alpha}}\right)^\gamma \int_D q(t-s, x, z) dz.
\end{eqnarray*}
Applying \eqref
{e:hkein1}, we have proved the lemma. \qed

\begin{lemma}\label{l:i3}
For $k \ge 0$
and $(t, x, y)\in (0, 1]\times D\times D$  we have
\begin{equation}\label{e:i3-1}
|p^k(t, x, y)| \le p^0(t,x,y) \left(\big(C_0^2MN^{\alpha, \gamma}_{\mu,
F_1}(t)\big)^k +k \|F_1\|_\infty C_0^2M \big(C_0^2MN^{\alpha, \gamma}_{\mu,
F_1}(t)\big)^{k-1}  \right).
\end{equation}
\end{lemma}

\pf
We use induction on
$k \ge 0$. The $k=0$ case is obvious.
Suppose that \eqref{e:i3-1} is true for $k-1\ge 0$.  Recall that
\begin{eqnarray*}
V_{x,y}={\{(z,w)\in D\times D:
|x-y|\ge 4(|y-w|
\wedge |x-z|) \}}, \quad
U_{x,y}=(D \times D) \setminus
V_{x,y}.
\end{eqnarray*}
Applying \eqref{e:i0}, \eqref{e:hke}, \eqref{e:M2} and \eqref{e:M4},
we have by our induction hypothesis
\begin{eqnarray*}
&&
|p^{k}(t, x, y)|
\,\le\,\int^t_0 \left(\int_{D}p^0(t-s, x,
z)
|p^{k-1}(s, z, y)|
|\mu|(dz) \right)ds \\
&&+ \int^t_0 \left(\int_{
U_{x,y}} p^0(t-s, x, z)\frac{c(z, w) |F_1(z,
w)|}{|z-w|^{d+\alpha}}
| p^{k-1}(s, w, y) |dzdw \right)ds\\
&&+\int^t_0 \left(\int_{
V_{x,y}} p^0(t-s, x, z)\frac{c(z, w)|F_1(z,
w)|}{|z-w|^{d+\alpha}}
| p^{k-1}(s, w, y) | dzdw \right)ds\\
&\le& \left(\big(C_0^2MN^{\alpha, \gamma}_{\mu, F_1}(t)\big)^{k-1}
+(k-1) \|F_1\|_\infty C_0^2M \big(C_0^2MN^{\alpha, \gamma}_{\mu, F_1}(t)\big)^{k-2} \right) \\
&&\qquad\times \int^t_0 \left(\int_{D}p^0(t-s, x, z)p^{0}(s, z,
y)|\mu|(dz) \right)  ds \\
&&+\left(\big(C_0^2MN^{\alpha, \gamma}_{\mu, F_1}(t)\big)^{k-1}
+(k-1)\|F_1\|_\infty C_0^2M \big(C_0^2MN^{\alpha,
\gamma}_{\mu, F_1}(t)\big)^{k-2} \right)  \\
&&\qquad\times  \int^t_0 \left(\int_{
U_{x,y}} p^0(t-s, x,
z)\frac{c(z, w)|F_1|(z,
w)}{|z-w|^{d+\alpha}}p^{0}(s, w, y) dzdw \right)ds     \\
&&+ \int^t_0 \left(\int_{
V_{x,y}} p^0(t-s, x, z)\frac{c(z, w)|F_1|(z,
w)}{|z-w|^{d+\alpha}}|p^{k-1}(s, w, y)| dzdw \right)ds\\
&\le& p^0(t,x,y) \left(\big(C_0^2MN^{\alpha, \gamma}_{\mu, F_1}(t)\big)^{k-1}
+(k-1) \|F_1\|_\infty C_0^2M \big(C_0^2MN^{\alpha, \gamma}_{\mu, F_1}(t)\big)^{k-2}
\right) C_0^2 M N^{\alpha, \gamma}_{\mu, F_1}(t)  \\
&&+C_0\frac{2^{d+\alpha}\|F_1\|_\infty}{|x-y|^{d+\alpha}} \int^t_0 \left(\int_{D\times D} p^0(t-s, x, z)
|p^{k-1} (s, w, y)|
dzdw \right)ds.
\end{eqnarray*}
Applying Lemma \ref{l:i2} and using \eqref{e:M1},
we get that if $|x-y|^\alpha \ge t$,
\begin{eqnarray*}
&&C_0\frac{2^{d+\alpha}\|F_1\|_\infty}{|x-y|^{d+\alpha}}\int^t_0
\left(\int_{D\times D} p^0(t-s, x, z)
|p^{k-1} (s, w, y)|
 dzdw \right)ds\\
&\le& \psi_\gamma (t,x,y) \frac{ t }{|x-y|^{d+\alpha}}
\|F_1\|_\infty
C_0^5\frac{\alpha}{\alpha-\gamma}
 2^{d+\alpha+\gamma/\alpha} M^{k-1}  \big(N^{\alpha, \gamma}_{\mu, F_1}(t)\big)^{k-1}\\
&\le&  p^0 (t,x,y)
\|F_1\|_\infty
C_0^6\frac{\alpha}{\alpha-\gamma}
 2^{d+\alpha+\gamma/\alpha} M^{k-1}  \big(N^{\alpha, \gamma}_{\mu, F_1}(t)\big)^{k-1}\\
&\le& p^0 (t,x,y) \|F_1\|_\infty C_0^2 M^{k}  \big(N^{\alpha,
\gamma}_{\mu, F_1}(t)\big)^{k-1} \\
&\le&  p^0 (t,x,y) \|F_1\|_\infty C_0^2 M  \big(C_0^2 M N^{\alpha,
\gamma}_{\mu, F_1}(t)\big)^{k-1}.
\end{eqnarray*}
Thus  \eqref{e:i3-1} is true for $k$ when  $|x-y|^\alpha \ge t$.

If $|x-y|^\alpha \le t$,
using \eqref{e:hke}, \eqref{e:i0}, \eqref{e:M2} and \eqref{e:M3},
we have by our induction hypothesis
\begin{eqnarray*}
&&
|p^{k}(t, x, y)|
\,\le\,\int^t_0 \left(\int_{D}p^0(t-s, x,
z)
|p^{k-1}(s, z, y)|
|\mu|(dz) \right)ds \\
&&+ \int^t_0 \left(\int_{
D \times D} p^0(t-s, x, z)\frac{c(z, w) |F_1(z,
w)|}{|z-w|^{d+\alpha}}
| p^{k-1}(s, w, y) |dzdw \right)ds\\
&\le& \left(\big(C_0^2MN^{\alpha, \gamma}_{\mu, F_1}(t)\big)^{k-1}
+(k-1) \|F_1\|_\infty C_0^2M \big(C_0^2MN^{\alpha, \gamma}_{\mu, F_1}(t)\big)^{k-2} \right)  \\
&&\qquad\times \int^t_0 \left(\int_{D}p^0(t-s, x, z)p^{0}(s, z,
y)|\mu|(dz) \right)  ds \\
&&+\left(\big(C_0^2MN^{\alpha, \gamma}_{\mu, F_1}(t)\big)^{k-1}
+(k-1)\|F_1\|_\infty C_0^2M \big(C_0^2MN^{\alpha,
\gamma}_{\mu, F_1}(t)\big)^{k-2} \right)   \\
&&\qquad\times  \int^t_0 \left(\int_{
D \times D} p^0(t-s, x,
z)\frac{c(z, w)|F_1|(z,
w)}{|z-w|^{d+\alpha}}p^{0}(s, w, y) dzdw \right)ds     \\
&\le& p^0(t,x,y) \left(\big(C_0^2MN^{\alpha, \gamma}_{\mu, F_1}(t)\big)^{k-1}
+(k-1) \|F_1\|_\infty C_0^2M \big(C_0^2MN^{\alpha, \gamma}_{\mu, F_1}(t)\big)^{k-2}
\right) C_0^2 M N^{\alpha, \gamma}_{\mu, F_1}(t).
\end{eqnarray*}
The proof is now complete.
\qed

Since $F_1 \in {\bf J}_{\alpha, \gamma}$, there is
$t_1 :=t_1(d, \alpha,
\gamma, C_0, M, N^{\alpha, \gamma}_{\mu, F_1}, \|F_1\|_\infty)\in (0, 1)$ so that
$$
N^{\alpha, \gamma}_{\mu, F_1}(t_1) \le
\big(3C_0^2M\big)^{-1} \wedge \big(9(C_0^2M)^2
\|F_1\|_\infty\big)^{-1}.
$$
It follows from Lemma \ref{l:i3} that
for every $(t, x,y) \in (0, t_1] \times D \times D$,
\begin{eqnarray}
&&\sum_{k=0}^\infty  |p^k(t, x, y)| \,=\,  p^0(t, x, y) +
\sum_{k=1}^\infty  |p^k(t, x, y)|\nonumber\\
&\le& p^0(t, x, y) +p^0(t,x,y) \left( \sum_{k=1}^\infty
(C_0^2MN^{\alpha, \gamma}_{\mu, F_1}(t))^k +\|F_1\|_\infty  C_0^2M
\sum_{k=1}^\infty  k
(C_0^2MN^{\alpha, \gamma}_{\mu, F_1}(t))^{k-1} \right) \nonumber\\
&\le& p^0(t, x, y) +p^0(t,x,y) \left( \frac12
+\frac94\|F_1\|_\infty  C_0^2M  \right) \nonumber \\
&\le&   p^0(t,x,y) \left( \frac32
+\frac94\|F_1\|_\infty  C_0^2M  \right). \label{e:newq1}
\end{eqnarray}
Hence, Fubini's Theorem, \eqref{e:newq3} and \eqref{e:1.12}
yield \eqref{e:1.11} and \eqref{e:1.14}.
Thus we conclude from
\eqref{e:newq1} and \eqref{e:1.14} that

\begin{thm}\label{t:ub}
There exist $t_1:=t_1
(d, \alpha, \gamma, C_0, M, N^{\alpha, \gamma}_{\mu, F_1}, \|F_1\|_\infty)\in (0, 1)$
and  a positive constant $C_6:=C_6
(d, \alpha,
\gamma, C_0, M, \|F_1\|_\infty)$ such that
the Feynman-Kac semigroup $T^{\mu, F}_t$
corresponding to
$\mu$ and $F$  has a continuous
density $q_D(t, x, y)$ for $t \le t_1$ and
\begin{equation}\label{e:hke1}
q_D(t, x, y) \le C_6 \psi_\gamma(t, x, y)q(t,x,y)
\end{equation}
for every $(t,x,y) \in (0, t_1] \times D \times D$.
\end{thm}

For the lower bound estimate, we need to assume that
$F$ is a   function in ${\bf J}_{\alpha, \gamma}$.

\begin{thm}\label{t:lb-1}
Suppose that
$\mu \in {\bf K}_{\alpha, \gamma}$
 and $F$ is a   function in
${\bf J}_{\alpha, \gamma}$.
Then there exist constants $t_2:=t_2(\alpha, \gamma, C_0, M, N^{\alpha, \gamma}_{\mu, F}, \|F\|_\infty)\in (0, 1)$ and
$C_7:=C_7(\alpha, \gamma, C_0, M, N^{\alpha, \gamma}_{\mu, F}, \|F\|_\infty)>1$
such that
\begin{equation}\label{e:hke2}
C_7^{-1} \psi_\gamma(t, x, y)q(t,x,y) \leq q_D(t, x, y)
\leq
C_7 \psi_\gamma(t, x, y)q(t,x,y)
\end{equation}
for every $(t,x,y) \in (0, t_2] \times D \times D$.
\end{thm}

\pf
Since $F$ is a bounded function in ${\bf J}_{\alpha, \gamma}$,
so is $F_1:=e^F-1$ with
$|F_1 (x, y) | \leq e^{\| F\|_\infty} |F| (x, y)$
and  $N^{\alpha, \gamma}_{F_1}\leq e^{\| F\|_\infty} N^{\alpha, \gamma}_{F}$.
Thus the upper bound estimate in \eqref{e:hke2}
follows directly from Theorem \ref{t:ub}.
To establish the lower bound, we define
for $(t, x, y)\in (0, \infty)\times D\times D$,
\begin{eqnarray*}
{\tilde p}^1(t, x, y)&=&\int^t_0\left(\int_Dp^0(t-s, x,
z)p^0(s, z, y)|\mu|(dz) \right)ds\\
&&+  \int^t_0\left(\int_D\int_Dp^0(t-s, x, z)\frac{c(z, w)|F|(z,
w)}{|z-w|^{d+\alpha}}p^0(s, w, y)dzdw \right)ds.
\end{eqnarray*}
Then for any bounded Borel function $f$ on $D$ and any $(t, x)\in
(0, \infty)\times D$, we have
$$
\E_x \left[ A^{|\mu|, |F|}_tf(X_t) \right]
=\int_D
{\tilde p}^1(t, x, y)f(y)dy.
$$
Applying Lemma \ref{l:i3} with $|\mu|$ and $|F|$ in place of
$\mu$ and $F_1$, we have
$$
{\tilde p}^1(t, x, y)\,\le\,(C^2_0MN^{\alpha, \gamma}_{\mu, F}(1)+C^2_0M\|F\|_\infty) p^0(t, x, y)
\,=:\,(k/2)p^0(t, x, y)
$$
for all $(t, x , y)\in
(0,1]\times D\times D$. Hence we have for
all $(t, x , y)\in
(0,1]\times D\times D$,
\begin{equation}\label{e:lbd1}
p^0(t, x, y)- \frac1k
{\tilde p}^1(t, x, y)\ge \frac12p^0(t, x,
y).
\end{equation}
Using the elementary fact that
$$
1-A^{|\mu|, |F|}_t/k\le \exp\left(-A^{|\mu|, |F|}_t/k\right)\le \exp\left( A^{\mu, F}_t/k \right),
$$
we get that for any $B(x, r)\subset D$ and any $(t, y)\in
(0, 1]\times D$,
$$
\frac1{|B(x, r)|}\E_y\left[\big(1-A^{|\mu|, |F|}_t/k \big)
{\bf 1}_{B(x,
r)}(X_t) \right] \le \frac1{|B(x, r)|}
\E_y\left[\exp (A^{\mu, F}_t/k)
{\bf 1}_{B(x, r)}(X_t) \right].
$$
Thus, by \eqref{e:lbd1} and H\"{o}lder's inequality, we have
\begin{eqnarray*}
&&\frac12\frac1{|B(x, r)|}\E_y\left[
{\bf 1}_{B(x, r)}(X_t) \right]\,\le \, \frac1{|B(x, r)|}\E_y\left[\exp (A^{\mu, F}_t/k)
{\bf 1}_{B(x,
r)}(X_t) \right]\\
&&\le
\left(\frac1{|B(x, r)|}\E_y\left[\exp (A^{\mu, F}_t)
{\bf 1}_{B(x, r)}(X_t) \right] \right)^{1/k}
\left(\frac1{|B(x, r)|} \E_y\left[
{\bf 1}_{B(x, r)}(X_t) \right]\right)^{1-1/k}.
\end{eqnarray*}
Therefore
$$
\frac1{2^k}\frac1{|B(x, r)|} \E_y\left[
{\bf 1}_{B(x, r)}(X_t) \right]
\le \frac1{|B(x, r)|}
 \E_y\left[\exp (A^{\mu, F}_t )
{\bf 1}_{B(x, r)}(X_t)
\right].
$$
We conclude by sending $r\downarrow 0$ that for every
$(t, x , y)\in
(0,1]\times D\times D$,
$2^{-k}p^0(t, x, y)\le q_D(t, x, y)$.
\qed

Combining the two theorems above with the semigroup property,
we immediately get the main result of this paper,
Theorem \ref{t:main}.

\section{Applications}\label{s:4}

In this section,
we will apply our main result to
(reflected) symmetric stable-like processes, killed symmetric
$\alpha$-stable processes, censored $\alpha$-stable processes  and stable processes with
drifts.
We first record the following two facts.

Suppose that $d\geq 2$ and $\alpha \in (0, 2)$.
A signed measure $\mu$  on $\bR^d$ is said to be in Kato class
$\bK_{d, \alpha}$
  $$\lim_{r\to 0} \sup_{x\in \bR^d} \int_{B(x, r)} \frac1{|x-y|^{d-\alpha}} |\mu| (dy) =0.
$$
A function $
g$ on $\bR^d$ is said to be in
$\bK_{d, \alpha}$ if $
g(x) dx\in \bK_{d, \alpha}$.

\begin{prop}\label{P:4.1} Suppose that
$d\geq 2$ and $\alpha \in (0, 2)$.

\begin{description}
\item{\rm (i)}
Let $D$ be an arbitrary Borel subset of $\bR^d$.
$\mu \in {\bf K}_{\alpha, 0}$
 if and only if ${\bf 1}_D  \mu  \in
\bK_{d, \alpha}$. Hence for every $\mu \in \bK_{d, \alpha}$, $\mu |_D
\in  {\bf K}_{\alpha, \gamma}$ for every $\gamma \geq 0$.
In particular,
$L^\infty (D; dx) \subset {\bf K}_{\alpha, \gamma}$
and $L^p(D; dx) \subset {\bf K}_{\alpha, \gamma}$ for every
$p>d/\alpha$ and $\gamma \geq 0$.

\item{\rm (ii)} Suppose that $D$ is a bounded Lipschitz open set in $\bR^d$ and $\gamma \in (0, \alpha)$.
Let $
g$ be a function defined on $D$.
 If there exist constants $c>0$, $\beta\in (0, \, \gamma + (\alpha -\gamma)/d)$ and a compact subset $K$
of $D$ such that $
{\bf 1}_K(x)
g(x)  \in \bK_{d, \alpha}$ and
$$
|
g(x)|\le
c \delta_D (x)^{-\beta}  \quad \hbox{for } x\in D\setminus K,
$$
then $
g \in {\bf K}_{\alpha, \gamma}$.
\end{description}
\end{prop}

\pf (i)
By the proof of \cite[Theorem 2]{Zh}, we have that
$\mu\in \bK_{d, \alpha}$ if and only if
$$ \lim_{t\to 0} \sup_{x\in \bR^d} \int_0^t \int_{\bR^d}
 q(s, x, y) \mu (dy) ds =0.
$$
This implies that $\mu \in {\bf K}_{\alpha, 0}$
if and only if ${\bf 1}_D \mu \in \bK_{d, \alpha}$.
In particular we have for every $\mu \in \bK_{d, \alpha}$,
$\mu |_D \in   {\bf K}_{\alpha, \gamma}$
for every $\gamma \ge 0$. Clearly  $L^\infty (D; dx)\subset \bK_{d, \alpha}$.
Using H\"older's inequality, one concludes that
 $L^p(\bR^d; dx) \subset \bK_{d, \alpha}$ for every
$p>d/\alpha$.

(ii) Let $
g$ be a function defined on $D$ such that there exist constants $
c_1>0$, $\beta\in (0, \, \gamma + (\alpha -\gamma)/d)$ and a compact subset $K$
of $D$ so that $
{\bf 1}_K(x)
g(x)  \in \bK_{d, \alpha}$ and
$|
g(x)|\le
c_1\delta_D (x)^{-\beta}$ for $x\in D\setminus K$.
In view of (i), it suffices to show that $
{\bf 1}_{D\setminus K}
g\in
{\bf K}_{\alpha, \gamma}$. Note that
\begin{eqnarray}
&& \sup_{x\in D}\int^t_0\int_{D\setminus K}\left(1\wedge
\frac{\delta_D(y)}{s^{1/\alpha}}\right)^\gamma
q(s, x, y)|
g(y)| dy ds \nn \\
&\leq & c_1 \sup_{x\in D}\int^t_0\int_{D\setminus K}\left(1\wedge
\frac{\delta_D(y)}{s^{1/\alpha}}\right)^\gamma \delta_D (y)^{-\beta}
q(s, x, y) dy ds  \nn \\
&\le
& c_1 \sup_{x\in D}
\int_{D\setminus K} \left( \int_0^{\delta_D(y)^\alpha \wedge t}
\left( s^{-d/\alpha} \wedge \frac{s}{|x-y|^{d+\alpha}}\right) ds \right) \delta_D(y)^{-\beta} dy \nn \\
&&+
 c_1 \sup_{x\in D}
\int_{D\setminus K} \left( \int_{\delta_D(y)^\alpha \wedge t}^t
s^{-\gamma /\alpha}
\left( s^{-d/\alpha} \wedge \frac{s}{|x-y|^{d+\alpha}}\right) ds \right) \delta_D(y)^{\gamma -\beta} dy
\,
=:\, I+II. \label
{e:4.2-1}
\end{eqnarray}
Here
\begin{eqnarray}
I &\leq & c_1 \sup_{x\in D}
\left(\int_D \int_0^{\delta_D (y)^\alpha \wedge |x-y|^\alpha
\wedge t} \frac{s}{|x-y|^{d+\alpha}} ds \,\delta_D(y)^{-\beta} dy \right.
 \nn \\
&& +
\left. \int_D \int_{\delta_D (y)^\alpha \wedge |x-y|^\alpha
\wedge t}^{\delta_D(y)^\alpha \wedge t} s^{-d/\alpha} ds \, \delta_D(y)^{-\beta} dy  \right)\nn \\
&\leq& c_2 \sup_{x\in D} \int_D \left(
  \frac{(\delta_D (y)^\alpha \wedge |x-y|^\alpha
\wedge t)^2 \delta_D (y)^{-\beta}}{|x-y|^{d+\alpha}}
  +
{\bf 1}_{\{ |x-y|< \delta_D (y) \wedge t^{1/\alpha}\}} \,  \frac{ \delta_D (y)^{-\beta}}{ |x-y|^{d-\alpha}}  \right) dy  \nn \\
&\leq &  c_2 \sup_{x\in D}  \int_D
\left( \frac{( |x-y|
\wedge t^{1/\alpha})^{2\alpha-\beta}} {|x-y|^d}
+{\bf 1}_{\{ |x-y|< \delta_D (y) \wedge t^{1/\alpha}\}}
\frac{1}{ |x-y|^{d-\alpha+\beta}} \right) dy \nn \\
&\leq &
2c_2t^{(\alpha-\beta)/(2\alpha)} \sup_{x\in D} \int_D \frac{1}{|x-y|^{d-(\alpha-\beta)/2}} dy
= c_3\, t^{(\alpha-\beta)/(2\alpha)},
\label{e:4.3}
\end{eqnarray}
while
\begin{eqnarray}
II
&\leq & c_1 \sup_{x\in D}
\int_{D} \left( \int_{\delta_D(y)^\alpha \wedge t}^t
{\bf 1}_{\{s<|x-y|^\alpha\}} \frac{s^{1-\gamma /\alpha}}{|x-y|^{d+\alpha}}
 ds \right) \delta_D(y)^{\gamma -\beta} dy  \nn \\
 && + c_1 \sup_{x\in D}
\int_{D} \left( \int_{\delta_D(y)^\alpha \wedge t}^t
{\bf 1}_{\{s\geq |x-y|^\alpha\}} \, s^{-(d+\gamma )/\alpha}
 ds \right) \delta_D(y)^{\gamma -\beta} dy  \nn \\
&\leq & c_4 \sup_{x\in D} \int_{D}
\frac{\left(|x-y|\wedge t^{1/\alpha}\right)^{2\alpha-\gamma}}{|x-y|^{d+\alpha}} \,
{\bf 1}_{\{ \delta_D(y)< |x-y| \wedge t^{1/\alpha}\}} \,
 \delta_D(y)^{\gamma -\beta} dy \nn\\
&& + c_4 \sup_{x\in D}   \int_D  |x-y|^{\alpha -d -\gamma}
{\bf 1}_{\{  |x-y| \le t^{1/\alpha}\}}
\delta_D(y)^{\gamma -\beta} dy \nn\\
&\leq & c_4 t^{\delta/\alpha}
\sup_{x\in D} \int_D \frac{1}{|x-y|^{d-\alpha +\eps} \delta_D (y)^{\beta -\gamma} } \, dy,  \label{e:4.4}
\end{eqnarray}
where $\delta := (\alpha -\gamma -d(\beta -\gamma))/2>0$
and $\eps := (\alpha +\gamma -d(\beta -\gamma))/2>0$.
Note that $\eps + \delta= \alpha-d(\beta-\gamma)$ and
$\eps-\delta =\gamma$.
Let $p=d/(d-\alpha+\eps+\delta/2)$ and $q= d/(\alpha-(\eps + \delta/2))$ so that
$1/p+1/q=1$. Since $D$ is a bounded Lipschitz open set,
$p(d-\alpha+\eps)<d$ and $q(\beta-\gamma)<1$,
we have by Young's inequality,
$$
\sup_{x\in D} \int_D \frac{1}{|x-y|^{d-\alpha +\eps}
 \delta_D (y)^{\beta -\gamma} } \, dy
  \leq  \sup_{x\in D} \int_D \left(\frac1{p} \frac1{|x-y|^{p(d-\alpha+\eps)}}
+ \frac1{q} \frac{1}{\delta_D(y)^{q(\beta-\gamma)}} \right) dy
<\infty.
$$
This together with \eqref{e:4.2-1}
--\eqref{e:4.4} implies that
$\lim_{t\to 0} N^{\alpha, \gamma}_{
{\bf 1}_{D\setminus K} g(x)}(t)=0$;
that is, ${\bf 1}_{D\setminus K}
g \in {\bf K}_{\alpha, \gamma}$.
This completes the proof of the proposition.
\qed

\begin{prop}\label{p:newqw1}
Suppose $\gamma  \in [0, \alpha \wedge d)$
and
$|F|(z,w) \le A(|z-w|^\beta \wedge 1)$ for some
$A>0$ and $\beta>\alpha$. Then there exists $
C_8=
C_8 (\beta, d, \alpha, \gamma)>0$
such that for every arbitrary Borel subset $D$ of $\bR^d$,
\begin{equation}\label{e:works for F}
N^{\alpha, \gamma}_F(t) \,\le\,
C_8\, A t.
\end{equation}
This in particular implies that $F\in {\bf J}_{\alpha, \gamma}$.
\end{prop}

\pf By \eqref{e:ineq1}, we have that
\begin{align*}
&\int_0^t\int_{\bR^d \times \bR^d}q(s, y, w )\left(1
 +
\frac{ |z-w| \wedge t^{1/\alpha}}{|y-w|} \right)^\gamma
\frac{|F|(z,w) + |F|(w,z)}{|z-w|^{d+\alpha}} dzdwds\\
\le &2 A\left(\int_{\bR^d}  (|z|^\beta \wedge 1)  |z|^{-d-\alpha} dz  \right)
\int_0^t\int_{\bR^d}q(s, y, w )\left(1
 +
\frac{ |z-w|
\wedge t^{1/\alpha}}{|y-w|} \right)^\gamma dwds\\
\le &   c_1A \left(\int_{B(0, 1)}  \frac{dz}{|z|^{d+\alpha-\beta}} +
\int_{B(0, 1)^c}  \frac{dz}{|z|^{d+\alpha}}    \right) \int_0^t
\left(1+ \int_{D}  2^{d+\alpha}  \frac{st^{\gamma/\alpha }}
{(s^{1/\alpha}+|y-w|)^{d+\alpha} |y-w|^\gamma}  \right) dwds\\
\le &  c_2A  \int_0^t  \left(1+ t^{\gamma /\alpha} s
\int^{\infty}_{0}  \frac{r^{d-1} }{r^{\gamma} (s^{1/\alpha}+r)^{d+\alpha}} dr \right) ds\\
\le &  c_2A t +c_3A
\left(\int^\infty_{0}  \frac{u^{d-1-\gamma}  }{
(1+u)^{d+\alpha}}  du \right)  t^{\gamma /\alpha}  \int_0^t  s^{-\gamma
/\alpha} ds \le c_4A t
\end{align*}
where the assumption $\gamma  \in [0, \alpha \wedge d)$ is used the last inequality.
This establishes \eqref{e:works for F}.
\qed

\subsection{Stable-like processes on closed $d$-sets}\label{S:4.1}

A Borel subset
$D$
in $\bR^d$ with $d\ge 1$ is said to be a  $d$-set if
 there exist constants $r_0>0$, $
\sC_2>
\sC_1>0$ so that
\begin{equation}\label{eqn:dset}
\sC_1\, r^d\le | B(x,r) \cap
 D
|\le
\sC_2 \, r^d \qquad \hbox{for all }~x\in
D \hbox{ and } 0<r\leq r_0,
\end{equation}
where for a Borel set $A\subset \bR^d$, we use $|A|$ to denote
its Lebesgue measure.
The notion of a $d$-set arises both in the theory of function spaces
and in fractal geometry. It is known that if $
D$ is a $d$-set, then
so is its Euclidean closure $\overline
D$.
Every uniformly Lipschitz open set in $\bR^d$ is a $d$-set,
so is its Euclidean
closure. It is easy to check that the classical von Koch snowflake
domain in $\bR^2$ is an open $2$-set. A $d$-set can have very rough
boundary since every $d$-set with a subset of zero Lebesgue
measure removed is still
a  $d$-set.

Suppose that $D$ is a closed $d$-set $D\subset \bR^d$ and
$c(x, y)$ is a symmetric function on
$D \times D$ that is bounded
between two strictly positive constants $
\sC_4>
\sC_3>0$, that is,
\begin{equation}\label{eqn:1.4}
\sC_3 \leq c(x, y) \leq
\sC_4 \qquad \hbox{for  a.e. } x, y \in D.
\end{equation}
For $\alpha\in (0, 2)$, we define
\begin{eqnarray}
\sF  &=& \left\{ u \in L^2(D;
dx): \, \int_{D\times D}
\frac{(u(x)-u(y))^2}
{|x-y|^{d+\alpha}} \, dx \, dy < \infty \right\}
 \label{e:DF1} \\
\sE(u,v)&=& \frac 12 \int_{D\times D} (u(x)-u(y))(v(x)-v(y)) \frac{
c(x, y) } {|x-y|^{d+\alpha}}\, dx \, dy, \quad u, v \in \sF.
\label{e:DF2}
\end{eqnarray}
It is easy to check that $(\sE, \sF)$ is a regular Dirichlet form on
$L^2(D, dx)$ and therefore there is an associated
symmetric Hunt process $X$ on $D$ starting from every point in $D$ except for an
exceptional set that has zero capacity. The process $X$ is called a
symmetric $\alpha$-stable-like process on $D$ in \cite{CK}.
When $c(x, y)$ is a constant function, $X$ is the reflected
$\alpha$-stable process
appeared in \cite{BBC}.
Note that when $D=
\bR^d$ and
$c(x, y)$ is a constant function, then $X$ is nothing but a
symmetric $\alpha$-stable process on $\bR^d$.

It follows as a special case from \cite[Theorem 1.1]{CK} that
  the symmetric stable-like process  $X$ on a closed $d$-set in $\bR^d$
  has a H\"older continuous
transition density function $p (t, x, y)$ with respect to the Lebesgue
measure on $D$
that satisfies the estimate \eqref{e:hke} with $\gamma=0$
and the comparison constant $C_0$ depending only on
$d$,  $\alpha$, $r_0$
and the constants $
\sC_k$, $k=1, \cdots, 4$ in \eqref{eqn:dset} and \eqref{eqn:1.4}.
In particular, this implies that the process $X$ can be refined
so it can start from every point in $D$.
Thus as a special case of Theorem \ref{t:main},
we have the following.

\begin{thm}\label{T:4.3} Suppose that $X$ is a   symmetric
$\alpha$-stable-like process on
a closed $d$-set $D$ in $\bR^d$.
 Assume $\mu\in {\bf K}_{\alpha,0}$ and $F\in {\bf J}_{\alpha, 0}$.
 Let $q $ be the density of the Feynman-Kac semigroup of $X$
corresponding to $A^{\mu, F}$. For any $T>0$, there exists a constant $
C_{9}>1$ such that for
all $(t, x, y)\in (0, T]\times D\times D$,
$$
C_{9}^{-1} q(t,x,y)\le q (t, x, y) \le
C_{9} q(t,x,y).
$$
\end{thm}

\begin{remark}\label{R:4.4} \rm
Let $n\geq 1$ be an integer and $d\in (0, n]$.
In general,  a Borel subset $D$
in $\bR^n$ is said to be a  $d$-set if
 there exist a measure $\mu$ and  constants $r_0>0$, $ \sC_2>
\sC_1>0$ so that
\begin{equation}\label{eqn:generaldset}
\sC_1\, r^d\le \mu ( B(x,r) \cap D)
\le \sC_2 \, r^d \qquad \hbox{for all }~x\in
D \hbox{ and } 0<r\leq r_0,
\end{equation}
It is established in \cite{CK} that for every $\alpha \in (0, 2)$,
a symmetric $\alpha$-stable-like process $X$
can always be constructed on any closed $d$-set $D$ in $\bR^n$
via the Dirichlet form $(\sE, \sF)$ on $L^2(D; \mu)$ defined by
\eqref{e:DF1}--\eqref{e:DF2} but with the
$d$-measure $\mu (dx)$
in place of the Lebesgue measure $dx$ there.
Moreover by \cite[Theorem 1.1]{CK}, the process $X$
has a jointly H\"older continuous transition density function
$p (t, x, y)$ with respect to the $d$-measure
$\mu$ on $D$ that satisfies the estimate \eqref{e:hke} with $\gamma=0$.
The proof of Theorem \ref{t:main} also works for such process $X$;
in other words, Theorem \ref{T:4.3} continues to hold for such kind
of symmetric stable-like processes. \qed
\end{remark}

\subsection{Killed symmetric $\alpha$-stable processes}\label{ss:stable}

A symmetric $\alpha$-stable process $X$ in $\bR^d$ is a L\'evy
process whose characteristic function is given by
$ \E_0 \left[ \exp(i \xi \cdot X_t ) \right]
          =e^{-t |\xi|^\alpha} .$
It is well-known that the process $X$ has a L\'evy
intensity function
$J(x, y)={\cal A}(d, -\alpha)|x-y|^{-(d+\alpha)},$
where
\begin{equation}\label{e:cal a}
{\cal A}(d, -\alpha)= \alpha 2^{-1+\alpha} \Gamma(\frac{d+\alpha}2)
\pi^{-d/2}(\Gamma(1-\frac{\alpha}2))^{-1}.
\end{equation}
Here $\Gamma$ is the Gamma function defined by $\Gamma(\lambda):=
\int^{\infty}_0 t^{\lambda-1} e^{-t}dt$ for every $\lambda > 0$.
Let $X^D$ be the killed symmetric $\alpha$-stable process $X^D$ in a $C^{1,1}$ open set $D$.
It follows from \cite{CKS} that
$X^D$ satisfies the assumption of Section \ref{s:1} with
$\gamma=\alpha/2$.
Thus as a special case of Theorem \ref{t:main},
we have the following.

\begin{thm}\label{t:ks} Suppose that $X$ is a killed symmetric
$\alpha$-stable process in
a $C^{1,1}$ open set $D$.
 Assume $\mu\in {\bf K}_{\alpha,
\alpha/2}$ and $F\in {\bf J}_{\alpha, \alpha/2}$. Let $q_D$ be the
density of the Feynman-Kac semigroup of $X$
corresponding to
 $A^{\mu,
F}$. For any $T>0$, there exists a constant $
C_{10}>1$ such that for
all $(t, x, y)\in (0, T]\times D\times D$,
$$
C_{10}^{-1} \psi_{\alpha/2}(t, x, y)q(t,x,y)\le q_D(t, x, y) \le
C_{10}
\psi_{\alpha/2}(t, x, y)q(t,x,y).
$$
\end{thm}

Let $X^{m}$ be a relativistic $\alpha$-stable
process in $\bR^d$ with mass $m>0$, i.e., $X^{m}$ is
a L\'evy process in $\bR^d$ with
$$
\E_0\left[\exp (i
\xi \cdot X^{m}_t)
\right]=\exp \left(t\left(
m^\alpha-(|\xi|^2+m^{2})^{\alpha/2} \right) \right).
$$
$X^{m}$ has a L\'evy
intensity function
$J^m(x,y)= {\cal A}(d, -\alpha) \varphi (m^{1/\alpha}|x-y|)|x-y|^{-d-\alpha}$
where
\begin{equation}\label{e:4.6}
 \varphi (r):= 2^{-(d+\alpha)} \, \Gamma \left(
\frac{d+\alpha}{2} \right)^{-1}\, \int_0^\infty s^{\frac{d+\alpha}{
2}-1} e^{-\frac{s}{ 4} -\frac{r^2}{ s} } \, ds,
\end{equation}
which is decreasing and is a smooth function of $r^2$ satisfying
 $\varphi(0) = 1$ and
 \begin{equation}\label{e:4.2}
\varphi (r) \asymp e^{-r}(1+r^{(d+\alpha-1)/2} ) \qquad \hbox{on } [0,
\infty)
\end{equation}
(see \cite[pp. 276-277]{CS03b} for details).

Let $X^{m, D}$ be a killed relativistic $\alpha$-stable process in a bounded $C^{1,1}$ open set.
Define
\begin{eqnarray*}
 K^m_t := \exp \left(\sum_{0<s \le t} \ln(\varphi(m^{1/\alpha}
 |X^D_s-X^D_{s-}|))
+m\,(t\wedge\tau_D)\right).
\end{eqnarray*}
Since
$\int_{\bR^d} J(x,y)-J^m(x,y)dy=m$ for all $x\in \bR^d$
(see \cite{Ry}),
it follows from \cite[p.279]{CS03b} that $X^{m, D}$ can be obtained
 from   the killed
symmetric $\alpha$-stable process $X^D$ in $D$ through
the non-local Feynman-Kac transform $K^m_t$.
That is, $\E_x \big[ f(X^{m, D}_t) \big]
:= \E_x \left[ K^m_t f(X^D_t) \right]$
By \eqref{e:4.6},
for any $M>0$, there exists a constant $c=c(d, \alpha, M,
\text{diam} (D))>0$ such that  for all
$m\in (0, M]$,
$|\ln(\varphi(m^{1/\alpha} |x-y|))| \leq c  (|x-y|^2 \wedge 1) $
 and so, by Proposition \ref{p:newqw1},
$ F_m(x, y) :=\ln(\varphi(m^{1/\alpha} |x-y|)) \in {\bf J}_{\alpha, \alpha/2}$.
The constant function $m$ is in ${\bf K}_{\alpha,
\alpha/2}$ and so
 $N^{\alpha, \alpha/2}_{ m, F_m}(t)$
goes to zero as $t$ goes to zero uniformly on $m \in (0,
M]$.
Thus, as an application of Theorem \ref{t:main},
we arrive at the following result, which is the bounded open
set case of a more general result
recently obtained in \cite{CKS2} by a different method.

\begin{thm}\label{t:krs}
Suppose that $D$ is
a bounded $C^{1, 1}$
open set in $\bR^d$. For any $m>0$, let $p^m_D$ be the transition
density of the killed relativistic $\alpha$-stable process with
weight $m$ in $D$. For any $
M>0$ and $T>0$, there exists a constant
$
C_{11}>1$ such that for all $m\in (0,
M]$ and $(t, x, y)\in (0,
T]\times D\times D$,
$$
C_{11}^{-1} \psi_{\alpha/2}(t, x, y)q(t,x,y)\le p^m_D(t, x, y) \le
C_{11}
\psi_{\alpha/2}(t, x, y)q(t,x,y).
$$
\end{thm}

\subsection{Censored stable  processes}

Fix an open set $D$ in $\bR^d$ with $d\geq 1$.
Recall that ${\cal A} (d, -\alpha)$ is the constant defined in \eqref{e:cal a}.
Define a
bilinear form $\sE $ on $C_c^\infty(D)$ by
 \bee\label{e:csdf}
 \sE (u, v):= \frac1{2}  \int_D \int_D (u(x)-u(y))(v(x)-v(y))
\frac{
{\cal A}
(d, -\alpha)} {|x-y|^{d+\alpha}}dxdy , \quad u, v \in
C_c^\infty(D).
 \eee
Using Fatou's lemma, it is easy to check that the bilinear form $ (
\sE  , C^\infty_c (D))$ is closable in $L^2 (D, dx)$. Let $\sF$ be
the closure of  $C^\infty_c(D)$  under the Hilbert inner product
$\sE _1:=\sE +(\,\cdot\, ,\, \cdot\, )_{L^2(D, dx)}.$ As noted in
\cite{BBC}, $(\sE , \sF )$ is Markovian and hence a regular
symmetric Dirichlet form on $L^2(D, dx)$, and therefore there is an
associated symmetric Hunt process $Y=\{Y_t, t\ge 0, \P_x, x\in D\}$
taking values in $D$ (cf. Theorem 3.1.1 of \cite{FOT}). The process
$Y$
is the censored $\alpha$-stable  process in $D$  that is
studied in \cite{BBC}.
By \eqref{e:csdf}, the jumping kernel
$J(x, y)$ of the censored
$\alpha$-stable process $Y$ is given by
$$
J(x, y)=\frac{
{\cal A}
(d, -\alpha)}{|x-y|^{d+\alpha}} \qquad \hbox{for } x, y \in
D.
$$

 As a particular case of a more general result established in \cite[Theorem 1.1]{CKS1}, when $\alpha \in (1, 2)$ and $D$ is a $C^{1, 1}$ open subset
  of $\bR^d$, the censored $\alpha$-stable process on $D$ satisfies the assumption of Section \ref{s:1} with $\gamma=\alpha-1$. Thus as a special case of Theorem
\ref{t:main}, we have the following:

\begin{thm}\label{t:cs} Suppose that
$\alpha \in (1, 2)$ and that
 $Y$ is a censored stable  process in a $C^{1,1}$ open set  $D$.
Assume $\mu\in {\bf
K}_{\alpha, \alpha-1}$ and $F\in {\bf J}_{\alpha, \alpha-1}$. Let
$q_D$ be the density of the Feynman-Kac semigroup of $Y$
corresponding to  $A^{\mu, F}$. For any $T>0$, there exists a constant
$C_{12}>1$ such that for all $(t, x, y)\in (0, T]\times D\times D$,
$$
C_{12}^{-1} \psi_{\alpha-1}(t, x, y)q(t,x,y)\le q_D(t, x, y) \le
C_{12}
\psi_{\alpha-1}(t, x, y)q(t,x,y).
$$
\end{thm}

Similar to \cite{BBC}, we can define a censored relativistic
$\alpha$-stable process in $D$.
Alternatively, with
$$
K_t := \exp \left(\sum_{0<s \le t} \ln(\varphi(m^{1/\alpha} (|Y_{s-}- Y_s|))+{\cal A} (d, -\alpha)
\int^t_0 \int_D\frac{1-\varphi(m^{1/\alpha} |Y_s-y|)}{|Y_s-y|^{\alpha+d}} dy ds\right),
$$
if $D$ is a bounded $C^{1,1}$ open set, a censored relativistic
stable process $Y^m$ can also  be obtained from the censored stable process $Y$ through the
Feynman-Kac transform $K_t$.
That is, $ \E_x [ f(Y^m_t)]=\E_x [ K_t f(Y_t) ]$
(see \cite{CKi, CS03b}).
By an argument similar to that of Subsection \ref{ss:stable}, one can see that
$ F_m:=\ln (\varphi(m^{1/\alpha} |x-y|))  \in {\bf J}_{\alpha, \alpha/2}$.
Moreover, since
$$
g_m(x):=\int_D (1-\varphi(m^{1/\alpha} |x-y|))|x-y|^{-\alpha-d} dy \le \int_{\bR^d} (1-\varphi(m^{1/\alpha} |x-y|))|x-y|^{-\alpha-d} dy =m, $$
$g_m \in {\bf K}_{\alpha,\alpha/2}$ and
$N^{\alpha, \alpha/2}_{g_m, F_m}(t)$
goes to zero as $t$ goes to zero uniformly on $m \in (0, M]$.
Thus as a particular case of Theorem \ref{t:cs},
we have the following.

\begin{thm}\label{t:crs}
Suppose that $\alpha \in (1, 2)$  and that $D$ is
 a bounded $C^{1, 1}$
open set in $\bR^d$. For any $m>0$, let $q^m_D$ be the transition
density of the censored relativistic $\alpha$-stable process with
weight $m$ in $D$. For any $
M>0$ and $T>0$, there exists a constant
$C_{13}>1$ such that for all $m\in (0,
M]$ and $(t, x, y)\in (0,
T]\times D\times D$,
$$
C_{13}^{-1} \psi_{\alpha-1}(t, x, y)q(t,x,y)\le q^m_D(t, x, y) \le
C_{13}
\psi_{\alpha-1}(t, x, y)q(t,x,y).
$$
\end{thm}

In fact, Theorems \ref{t:cs} and \ref{t:crs} are applicable
to certain class of censored stable-like processes whose Dirichlet
heat kernel estimates are given in \cite{CKS1}.

\subsection{Stable processes with drifts}

Let $\alpha \in (1, 2)$ and $d\ge 2$. In this subsection, we apply our main result to a non-symmetric process.

For $b=(b_1, \dots, b_d)$ with $b_i \in {\mathbb K}_{d, \alpha-1}$,
a Feller process $Z$ on $\bR^d$ with infinitesimal generator
$\sL^b:=\Delta^{\alpha/2}+ b(x)\cdot \nabla$
is constructed in \cite{BJ} through the fundamental solution of $\sL^b$.
Let $Z^D$ be the subprocess of $Z$ killed upon leaving $D$.
The following result is established in \cite{CKS4}.

\begin{thm}
If $\alpha \in (1, 2)$, $d\ge 2$ and $D$
 is a bounded $C^{1, 1}$  open set, then $Z^D$ has a jointly continuous transition density
function $p_D(t, x, y)$ that satisfies \eqref{e:hke} with
$\gamma = \alpha/2$.
\end{thm}
Thus as a special case of Theorem \ref{t:main},
we also have the following:

\begin{thm}\label{t:ksnew}
Suppose that $\alpha \in (1, 2)$, $d\ge 2$, that
$D$ is a bounded $C^{1, 1}$  open set
and that  $Z^D$ is the subprocess of $Z$
killed upon leaving $D$. Assume $\mu\in {\bf K}_{\alpha,
\alpha/2}$ and $F\in {\bf J}_{\alpha, \alpha/2}$. Let $q_D$ be the
density of the Feynman-Kac semigroup of
$Z^D$ corresponding to  $A^{\mu,
F}$. For any $T>0$, there exists a constant $
C_{14}>1$ such that for
all $(t, x, y)\in (0, T]\times D\times D$,
$$
C_{14}^{-1} \psi_{\alpha/2}(t, x, y)q(t,x,y)\le q_D(t, x, y) \le
C_{14}
\psi_{\alpha/2}(t, x, y)q(t,x,y).
$$
\end{thm}

\vskip 0.3truein

{\bf Zhen-Qing Chen}

Department of Mathematics, University of Washington, Seattle,
WA 98195, USA

E-mail: \texttt{zqchen@uw.edu}

\bigskip

{\bf Panki Kim}

Department of Mathematics,
Seoul National University,
Building 27,
1 Gwanak-ro, Gwanak-gu,
Seoul 151-747, Republic of Korea

E-mail: \texttt{pkim@snu.ac.kr}

\bigskip

{\bf Renming Song}

Department of Mathematics, University of Illinois, Urbana, IL 61801, USA

E-mail: \texttt{rsong@math.uiuc.edu}
\end{document}